\documentclass{amsart}

 \newtheorem{thm}{Theorem}[section]
 
 \newtheorem{lem}[thm]{Lemma}
 
 \theoremstyle{definition}
 \newtheorem{defn}[thm]{Definition}
 \theoremstyle{remark}
 \newtheorem{rem}[thm]{Remark}
 \theoremstyle{definition}
 \newtheorem{exercise}[thm]{Exercise}
\newtheorem{example}[thm]{Example}

 \DeclareMathOperator{\Sub}{Sub}
 \DeclareMathOperator{\Seg}{Seg}

 \newcommand{\CC}{\mathbb{C}}
 \newcommand{\PP}{\mathbb{P}}

\def\move-in{\parshape=1.75true in 5true in}

\usepackage{xypic, amssymb}

\usepackage{url}

\begin{document}

\title{Four Lectures on Secant Varieties}

\author[E. Carlini]{Enrico Carlini}
\address[E. Carlini]{School of Mathematical Sciences, Monash University, Melbourne, VIC, 3800, Australia}
\email{enrico.carlini@monash.edu}
\address[E. Carlini]{Department of Mathematical Sciences, Politecnico di Torino, 10129, Turin, Italy}
\email{enrico.carlini@polito.it}

\author[N. Grieve]{Nathan Grieve}
\address[N. Grieve]{Department of Mathematics and Statistics,
McGill University,
Montreal, QC, Canada}
\email{ngrieve@math.mcgill.ca}
\address[N. Grieve]{Department of Mathematics and Statistics, Queen's University,
Kingston, ON, Canada}
\email{nathangrieve@mast.queensu.ca}

\author[L. Oeding]{Luke Oeding}
\address[L. Oeding]{Department of Mathematics and Statistics, Auburn University, Auburn, AL}
\email{oeding@auburn.edu}
\address[L. Oeding]{Department of Mathematics, University of California, Berkeley, Berkeley, CA}
\email{oeding@math.berkeley.edu}

\maketitle


\begin{abstract}
This paper is based on the first author's lectures at the 2012 University of Regina Workshop ``Connections Between Algebra and Geometry''.  Its aim is to provide an introduction to the theory of higher secant varieties and their applications.  Several references and solved exercises are also included.
\end{abstract}



\section{Introduction}

Secant varieties have travelled a long way from 19th century geometry to nowadays where they are as popular as ever before. There are different reasons for this popularity, but they can summarized in one word: applications.   These applications are both pure and applied in nature.  Indeed, not only does the geometry of secant varieties play a role in the study projections of a curve, a surface or a threefold, but it also in locating a transmitting antenna \cite{Comon2_tensor}. 

In these lectures we introduce the reader to the study of (higher) secant varieties by providing the very basic definition and properties, and then moving the most direct applications. In this way we introduce the tools and techniques which are central for any further study in the topic. In the last lecture we present some more advanced material and provide pointers to some 
relevant literature. Several exercises are included and they are meant to be a way for the reader to familiarize himself or herself with the main ideas presented in the lectures.  So, have fun with secant varieties and their many applications!

The paper is structured as follows. In Section \ref{lectone} we provide the basic definition and properties of higher secant varieties. In particular, we introduce one of the basic result in the theory, namely Terracini's Lemma, and one main source of examples and problems, namely Veronese varieties. In Section \ref{lecttwo} we introduce Waring problems and we explore the connections with higher secant varieties of Veronese varieties. Specifically, we review the basic result by Alexander and Hirschowitz. In Section \ref{lectthree}, Apolarity Theory makes in its appearance with the Apolarity Lemma. We see how to use Hilbert functions and sets of points to investigate Waring problems and higher secant varieties to Veronese varieties. In Section \ref{lectfour} we give pointers to the literature giving reference to the topics we treated in the paper. We also provide a brief description and references for the many relevant topics which we were not able to include in the lectures because of time constraints.  Finally, we provide solutions to the exercises in Section \ref{soltns}.

\medskip
\noindent
\emph{Acknowledgements.} The second and third authors would like to thank the first author for his lectures at the workshop ``Connections Between Algebra and Geometry" held at the University of Regina in 2012.  The second author would like to thank the third author for tutoring these lectures and for helping to write solutions to the exercises.  All three authors would like to thank the organizers S. Copper, S. Sather-Wagstaff and D. Stanley for their efforts in organizing the workshop and for securing funding to cover the costs of the participants.  It is also our pleasure to acknowledge the lecture notes of Tony Geramita \cite{MR1381732} which have influenced our exposition here. The three authors received partial support by different sources: Enrico Carlini by GNSAGA of INDAM, Nathan Grieve by an Ontario Graduate Fellowship, and Luke Oeding  by NSF RTG Award \# DMS-0943745. Finally, all the authors thank the anonymous referee for the improvement to the paper produced by the referee's suggestions and remarks.



\section{Lecture One}\label{lectone}

In what follows, $X\subset\mathbb{P}^N$ will denote an
irreducible, reduced algebraic variety; we work over an
algebraically closed field of characteristic zero, which we assume
to be $\mathbb{C}$.

The topic of this lecture are {\it higher secant varieties} of $X$.

\begin{defn} The $s$-th higher secant variety of $X$ is
\[\sigma_s(X)=\overline{\bigcup_{P_1,\ldots,P_s\in X}\langle P_1,\ldots,P_s\rangle},\]
where the over bar denotes the Zariski closure.
\end{defn}
In words, $\sigma_s(X)$ is the closure of the union of $s$-secant
spaces to $X$. \footnote{Some authors use the notation $S(X)$ for the (first) secant variety of $X$, which corresponds to $\sigma_{2}(X)$, and $S_{k}(X)$ to denote the $k$th secant variety to $X$, which corresponds to $\sigma_{k+1}(X)$. We prefer to reference the number of points used rather than the dimension of their span because this is often more relevant for applications because of its connection to rank.}

\begin{example} If $X\subset\mathbb{P}^2$ is a curve and not a line
then $\sigma_2(X)=\mathbb{P}^2$, the same is true for
hypersurfaces which are not hyperplanes. But, if
$X\subset\mathbb{P}^3$ is a non-degenerate curve (\emph{i.e.} not
contained in a hyperplane), then $\sigma_2(X)$ can be, in
principle, either a surface or a threefold.
\end{example}

We note that the closure operation is in general necessary, but
there are cases in which it is not.

\begin{exercise}\label{exONE1} Show that the union of chords (secant lines) to a
plane conic is closed. However, the union of the chords of the
twisted cubic curve in $\mathbb{P}^3$ is not.
\end{exercise}

In general,  we have a sequence of inclusions
\[X=\sigma_1(X) \subseteq \sigma_2(X) \subseteq \ldots \subseteq \sigma_r(X)\subseteq \ldots \subseteq \mathbb{P}^N.\]
If $X$ is a linear space, then $\sigma_i(X)=X$ for all $i$ and all
of the elements of the sequence are equal.

\begin{rem} If $X=\sigma_2(X)$ then  $X$ is a linear space. To
see this consider a point $P\in X$ and the projection map
$\pi_P:\mathbb{P}^N\dashrightarrow\mathbb{P}^{N-1}$. Let
$X_1=\pi_P(X)$ and notice that $\dim X_1=\dim X-1$ and that
$\sigma_2(X_1)=X_1$. If $X_1$ is a linear space also $X$ is so and
we are done. Otherwise iterate the process constructing a sequence
of varieties $X_2,\ldots, X_m$ of decreasing dimension. The
process will end with $X_m$ equal to a point and then $X_{m-1}$ a
linear space. Thus $X_{m-2}$ is a linear space and so on up to the
original variety $X$.
\end{rem}

\begin{exercise}\label{exONE2} For $X\subset\mathbb{P}^N$, show that, if
$\sigma_i(X)=\sigma_{i+1}(X)\neq\mathbb{P}^N$, then $\sigma_i(X)$
is a linear space and hence $\sigma_j(X)=\sigma_i(X)$ for all
$j\geq i$.
\end{exercise}

Using this remark and Exercise \ref{exONE2}, we can refine our chain of
inclusions for $X$ a non degenerate variety (\emph{i.e.} not contained in
a hyperplane).
\begin{exercise}\label{exONE2+}  If $X \subseteq \mathbb{P}^N$ is non-degenerate, then there exists an $r \geq 1$ with the property that
\[X = \sigma_1(X) \subsetneq \sigma_2(X) \subsetneq \ldots \subsetneq \sigma_r(X) = \mathbb{P}^N.\]
In particular, all inclusions are strict and there is a higher
secant variety that coincides with the ambient space.
\end{exercise}

It is natural to ask: what is the smallest $r$ such that
$\sigma_r(X)=\mathbb{P}^N$? Or more generally: what is the value
of $\dim \sigma_i(X)$ for all $i$?

As a preliminary move in this direction, we notice that there is
an expected value for the dimension of any higher secant variety
of $X$ that arises just from the naive dimension count. That is to say, if the secant variety doesn't fill the ambient space, a point is obtained by choosing $s$ points from an $n$-dimensional variety and one point in the $\PP^{s-1}$ that they span.

\begin{defn}For $X\subset\mathbb{P}^N$, set $n=\dim X$. The {\it expected dimension} of
$\sigma_s(X)$ is
\[\mathrm{expdim}(\sigma_s(X))=\min\{sn+s-1,N\}.\]
\end{defn}
Notice also that the expected dimension is also the maximum dimension of the secant variety.  Moreover, if the secant line variety  $\sigma_{2}(X)$ does not fill the ambient $\PP^{N}$, then $X$ can be isomorphically projected into a $\PP^{N-1}$.  This interest in minimal codimension embeddings is one reason secants were classically studied.

\begin{exercise}\label{exONE2++}
Let $X \subseteq \PP^n$ be a curve.  Prove that $\sigma_2(X)$ has dimension $3$ unless $X$ is contained in a plane.  (This is why every curve is isomorphic to a space curve but only birational to a plane curve.)
\end{exercise}

There are cases in which
$\mathrm{expdim}(\sigma_i(X))\neq\dim(\sigma_i(X))$ and these
motivate the following

\begin{defn} If
$\mathrm{expdim}(\sigma_i(X))\neq\dim(\sigma_i(X))$ then $X$ is
said to be {\it $i$-defective} or simply {\it defective}.
\end{defn}

\begin{rem} Notice that $\dim(\sigma_{i+1}(X))\leq
\dim(\sigma_{i}(X))+n+1$, where $n=\dim X$. This means that if
$\sigma_i(X)\neq\mathbb{P}^N$ and $X$ is $i$-defective, then $X$
is $j$-defective for $j\leq i$.
\end{rem}

Let's now see the most celebrated example of a defective variety,
the Veronese surface in $\mathbb{P}^5$.
\begin{example} Consider the polynomial ring $S=\mathbb{C}[x,y,z]$
and its homogeneous pieces $S_d$. The Veronese map $\nu_2$ is the morphism
\[\nu_2 \colon \mathbb{P}(S_1)\longrightarrow\mathbb{P}(S_2)  \text{ defined by }
[L]\mapsto [L^2] \text{.}\]
In coordinates this map can be described in terms of the standard
monomial basis $\langle x,y,z \rangle $ for $S_1$ and the standard monomial basis
$ \langle x^2,2xy,2xz,y^2,2yz,z^2\rangle $ for $S_2$.
Thus the Veronese map can be written as the map
\[\nu_2\colon \mathbb{P}^2\longrightarrow\mathbb{P}^5 \text{ defined by }
[a:b:c]\mapsto [a^2:ab:ac:b^2:bc:c^2] \text{.}\]
The Veronese surface is then defined as the
image of this map, \emph{i.e.} the Veronese surface is
$X = \nu_2(\mathbb{P}^2) \subseteq \PP^5$.

We now want to study the higher secant varieties of the Veronese surface $X$.  In
particular we ask: is $\dim\sigma_2(X) = \mathrm{expdim}
(\sigma_2(X)) = 5$? In other words, does $\sigma_2(X) $ equal $\mathbb{P}^5$?

To answer this question, it is useful to notice that elements in $S_2$ are quadratic forms
and, hence, are uniquely determined by $3\times 3$ symmetric
matrices. In particular, $P\in\mathbb{P}^5$ can be seen as $P=[Q]$
where $Q$ is a $3\times 3$ symmetric matrix. If $P\in X$ then $Q$
also has  rank equal one. Thus we have,
\[\sigma_2(X)=\overline{\bigcup_{P_1,P_2}\langle
P_1,P_2\rangle}\]
\[=\overline{\lbrace [Q_1+Q_2] : Q_i\mbox{ is a $3\times 3$ symmetric matrix and } \mathrm{rk}(Q_i)=1 \rbrace}
\subseteq H \text{,} \] where $H$ is the projective variety determined by the set of $3\times 3$ symmetric matrices of rank at most two.
Clearly $H$ is the hypersurface defined by the vanishing of the
determinant of the general $3\times 3$ symmetric matrix and hence
$X$ is $2$-defective.
\end{example}

\begin{exercise}\label{exONE2+++} Let $M$ be an $n\times n$ symmetric matrix of rank $r$.  Prove that $M$ is a sum of $r$ symmetric matrices of rank $1$.
\end{exercise}

\begin{exercise}\label{exONE3} Show that $H=\sigma_2(X)$.
\end{exercise}

\begin{exercise}\label{exONE4} Repeat the same argument for
$X=\nu_2(\mathbb{P}^3)$. Is $X$ $2$-defective?
\end{exercise}

In order to deal with the problem of studying the dimension of the
higher secant varieties of $X$ we need to introduce a celebrated
tool, namely Terracini's Lemma, see \cite{Terracini}.

\begin{lem}[Terracini's Lemma] Let $P_1,\ldots,P_s\in X$ be
general points and $P\in\langle
P_1,\ldots,P_s\rangle\subset\sigma_s(X)$ be a general point. Then
the tangent space to $\sigma_s(X)$ in $P$ is
\[T_P(\sigma_s(X))=\langle T_{P_1}(\sigma_s(X)),\ldots, T_{P_s}(\sigma_s(X))\rangle.\]
\end{lem}

\begin{rem} To get a (affine) geometric idea of why Terracini's
Lemma holds, we consider an affine curve $\gamma(t)$. A general
point on $P\in\sigma_2(\gamma)$ is described as
$\gamma(s_0)+\lambda_0[\gamma(t_0)-\gamma(s_0)]$. A neighborhood
of $P$ is then described as
$$\gamma(s)+\lambda[\gamma(t)-\gamma(s)] \text{.}$$ Hence the tangent space
$T_P(\sigma_s(\gamma))$ is spanned by
\[\gamma'(s_0)-\lambda_0\gamma'(s_0),\lambda_0\gamma'(t_0),\gamma(t_0)-\gamma(s_0),\]
and this is the affine span of the affine tangent spaces
$\{\gamma(s_0)+\alpha\gamma'(s_0): \alpha\in\mathbb{R}\}$ and
$\{\gamma(t_0)+\beta\gamma'(t_0): \beta\in\mathbb{R}\}$.
\end{rem}

As a first application of Terracini's Lemma, we consider the
twisted cubic curve.

\begin{example} Let $X$ be the twisted cubic curve in
$\mathbb{P}^3$, \emph{i.e.} $X=\nu_{3}(\mathbb{P}^1)$ where $\nu_{3}$ is the map
\[\nu_{3}\colon \mathbb{P}^1\longrightarrow\mathbb{P}^3 \text{ defined by } [s:t]\mapsto [s^3:s^2t:st^2:t^3].\]

We want to compute $\dim\sigma_2(X)=\dim T_P(\sigma_2(X))$ at a
generic point $P$. Using Terracini's Lemma it is enough to choose
generic points $P_1,P_2\in X$ and to study the linear span
\[\langle T_{P_1}(X),T_{P_2}(X)\rangle.\]
In particular, $\sigma_2(X)=\mathbb{P}^3$ if and only if the lines
$T_{P_1}(X)$ and $T_{P_2}(X)$ do not intersect, that is, if and
only if there does not exist a hyperplane containing both lines.

If $H\subset \mathbb{P}^3$ is a hyperplane, then the points of $H\cap X$
are determined by finding the roots of the degree three
homogeneous polynomial $g(s,t)$ defining
$\nu_{3}^{-1}(H)\subset\mathbb{P}^1$. If $H\supset T_{P_1}(X)$ then
$g$ has a double root. However, the homogeneous polynomial is smooth and thus, in the general case, no hyperplane exists containing both tangent lines.

In conclusion, $\sigma_2(X)=\mathbb{P}^3$.

\end{example}

\begin{exercise}\label{exONE5} Prove that if $H\supset T_{P}(X)$, then the
polynomial defining $\nu_{3}^{-1}(H)$ has a double root.
\end{exercise}

We now introduce the Veronese variety in general.

\begin{defn}
 Consider the polynomial ring $S=\mathbb{C}[x_0,\ldots,x_n]$
and its homogeneous pieces $S_d$. The {\it $d$-th Veronese map}
$\nu_d$ is the morphism
\[\nu_d \colon \mathbb{P}(S_1)\longrightarrow\mathbb{P}(S_d) \text{ defined by }
 [L]\mapsto [L^d].\]
In coordinates, using suitable monomial bases for $S_1$ and $S_d$, $\nu_d$ is the morphism
\[\nu_d \colon \mathbb{P}^n\longrightarrow\mathbb{P}^N \text{ defined by } \]
\[[a_0: \ldots: a_n]\mapsto [M_0(a_0,\dots, a_n) : \ldots : M_{N}(a_0,\dots, a_n)]\]
where $N={n+d\choose d}-1$ and where $M_0,\dots, M_N$ are monomials which form a basis for $S_d$.

We call $\nu_d(\mathbb{P}^n)$ a {\it Veronese variety}.
\end{defn}

\begin{example} A relevant family of Veronese varieties are the {\it rational normal
curves} which are Veronese varieties of dimension one, \emph{i.e.} $n=1$.
In this situation $S=\mathbb{C}[x_0,x_1]$ and $S_d$ is the vector
space of degree $d$ binary forms. The rational normal curve
$\nu_d(\mathbb{P}(S_1))\subseteq\mathbb{P} (S_d)$ is represented by
$d$-th powers of binary linear forms.
\end{example}

\begin{example} The rational normal curve  $X=\nu_2(\mathbb{P}^1)\subset\mathbb{P}^2$ is an
irreducible conic. It is easy to see that
$\sigma_2(X)=\mathbb{P}^2=\mathbb{P}(S_2)$. This equality can also
be explained by saying that any binary quadratic form $Q$ is the
sum of two squares of linear forms, \emph{i.e.} $Q=L^2+M^2$.
\end{example}

\begin{exercise}\label{exONE6} Consider the rational normal curve in
$\mathbb{P}^3$, \emph{i.e.} the twisted cubic curve $X=\nu_3(\mathbb{P}
(S_1))\subset\mathbb{P} (S_3)$. We know that $\sigma_2(X)$ fills up
all the space. Can we write any binary cubic as the sum of two
cubes of linear forms? Try $x_0x_1^2$.
\end{exercise}

\begin{exercise}\label{exONE7} We described the Veronese variety $X=\nu_d(\mathbb{P}^n)$ in parametric form
by means of the relation: $[F]\in X$ if and only if $F=L^d$. Use
this description and standard differential geometry to compute
$T_{[L^d]}(X)$ (describe this as a vector space of homogeneous
polynomials). This can be used to apply Terracini's Lemma, for
example, to the twisted cubic curve.
\end{exercise}


\section{Lecture Two}\label{lecttwo}

In the last lecture we spoke about higher secant varieties in
general. Now we focus on the special case of Veronese varieties.
Throughout this lecture we will consider the polynomial ring
$S=\mathbb{C}[x_0,\ldots, x_n]$.

An explicit description of the tangent space to a Veronese variety
will be useful, so we give it here.

\begin{rem} Let $X=\nu_d(\mathbb{P}^n)$ and consider
$P=[L^d]\in X$ where $L\in S_1$ is a linear form. Then
\[T_P(X)=\langle [L^{d-1}M] : M\in S_1\rangle.\]
\end{rem}

We can use this to revisit the Veronese surface example.

\begin{example} Consider the Veronese surface
$X=\nu_2(\mathbb{P}^2)\subset\mathbb{P}^5$. To compute
$\dim (\sigma_2(X))$ we use Terracini's Lemma. Hence we choose two
general points $P = [L^2],Q = [N^2]\in X$ and we consider the linear
span of their tangent spaces
\[T=\langle T_P(X),T_Q(X)\rangle.\]
By applying Grassmann's formula, and noticing that $T_P(X)\cap
T_Q(X)=[LN]$ we get $\dim T= 3 + 3 -1 -1=4$ and hence
$\sigma_2(X)$ is a hypersurface.
\end{example}

The study of higher secant varieties of Veronese varieties is
strongly connected with a problem in polynomial algebra: the {\it
Waring problem for forms}, \emph{i.e.} for homogeneous polynomials, see \cite{MR1381732}. We
begin by introducing the notion of Waring rank.

\begin{defn} Let $F\in S$ be a degree $d$ form. The {\it Waring rank of
$F$ } is denoted $\mathrm{rk}(F)$ and is defined to be the {\bf minimum} $s$
such that we have
\[F = L_1^d+\ldots +L_s^d\]
for some linear forms $L_i\in S_1$.
\end{defn}

\begin{rem} It is clear that $\mathrm{rk}(L^d)=1$ if $L$ is a linear
form.  However in general, if $L$ and $N$ are linear forms, $\mathrm{rk}(L^d+N^d)\leq 2$.  It is $1$ if $L$ and
$N$ are proportional and $2$ otherwise. For more than two factors
the computation of the Waring rank for a sum of powers of linear
form is not trivial.
\end{rem}

We can now state the {\it Waring problem for forms}, which
actually comes in two fashions. The {\bf big} Waring problem asks
for the computation of
\[g(n,d)\]
the minimal integer such that
\[\mathrm{rk}(F)\leq g(n,d)\]
for a {\bf generic} element $F \in S_d$, \emph{i.e.} for a generic
degree $d$ form in $n+1$ variables. The {\bf little} Waring
problem is more ambitious and asks us to determine the smallest
integer
\[G(n,d)\]
such that
\[\mathrm{rk}(F)\leq G(n,d)\] {\bf for any}
 $F \in S_d$.

\begin{rem} To understand the difference between
the big and the little Waring problem we can refer to a
probabilistic description. Pick a random element $F\in S_d$, then
$\mathrm{rk}(F) \leq G(n,d)$ and with {\it probability one}
$\mathrm{rk}(F) = g(n,d)$ (actually equality holds). However, if
the choice of $F$ is unlucky, it could be that $\mathrm{rk}(F)
> g(n,d)$. Note also that these notions are field dependent, see \cite{ComonOtt_Binary} for example.
\end{rem}

\begin{rem} To make the notion of a generic element precise we use
topology. Specifically, the big Waring problem asks us to bound the Waring rank
for all elements belonging to a non-empty Zariski open subset of
$\mathbb{P} S_d$; since non-empty Zariski open subsets are dense this also
explains the probabilistic interpretation.
\end{rem}

The big Waring problem has a nice geometric interpretation using
Veronese varieties --- this interpretation allows for a complete
solution to the problem. Also the little Waring problem has a
geometric aspect but this problem, in its full generality, is
still unsolved.

\begin{rem} As the Veronese variety $X=\nu_d(\mathbb{P}^n)\subset\mathbb{P}^N$
parameterizes pure powers in $S_d$, it is clear that $g(n,d)$ is
the smallest $s$ such that $\sigma_s(X)=\mathbb{P}^N$. Thus
solving the big Waring problem is equivalent to finding the
smallest $s$ such that secant variety $\sigma_s(X)$ fills up $\PP^N$.  On the other hand, as taking
the Zariski closure of the set $\bigcup_{P_1,\ldots,P_s\in X}\langle P_1,\ldots,P_s\rangle$ is involved in defining $\sigma_s(X)$, this is not equivalent to solving the little Waring problem.
\end{rem}

\begin{rem} To solve the little Waring problem one has to find the
smallest $s$ such that every element $[F]\in\mathbb{P}S_d$
lies on the span of some collection of $s$ points of $X$.
\end{rem}

Let's consider two examples to better understand the difference
between the two problems.

\begin{example} Let $X=\nu_2(\mathbb{P}^1)\subset\mathbb{P}^2$ be the rational normal
curve in $\mathbb{P}^2$, \emph{i.e.} a non-degenerate conic. We know that $\sigma_2(X) = \mathbb{P}^2$
and hence $g(n=1,d=2)=2$. But we also know that each point of
$\mathbb{P}^2$ lies on the span of two distinct points of $X$ --- every $2\times 2$ symmetric matrix is the sum of two rank-one symmetric matrices ---
thus $G(n=1,d=2)=2$. In particular this means that the Waring rank
of a binary quadratic form is always at most two.
\end{example}

\begin{example} Let $X=\nu_3(\mathbb{P}^1)\subset\mathbb{P}^3$ be the rational normal
curve. Again, we know that $\sigma_2(X)=\mathbb{P}^3$ and hence
$g(n=1,d=3)=2$. However, there are degree three binary forms $F$
such that $\mathrm{rk}(F)=3$, and actually $G(n=1,d=3)=3$.

To
understand which the bad forms are, consider the projection map
$\pi_P$ from any point $P=[F]\in\mathbb{P}^3$. Clearly, if
$P\not\in X$, $\pi_p(X)$ is a degree $3$ rational plane curve.
Hence, it is singular, and being irreducible, only two
possibilities arise. If the singularity is a node, then $P=[F]$
lies on a chord of $X$, and thus $F=L^3+N^3$. But, if the
singularity is a cusp, this is no longer true as $P$ lies on a
tangent line to $X$ and not on a chord. Thus, the bad binary
cubics lie on tangent lines to the twisted cubic curve. In other
words, the bad binary cubics are of the form $L^2N$.
\end{example}

\begin{exercise}\label{exTWO1} For binary forms, we can stratify $\mathbb{P} S_2$ using the Waring
rank: rank one elements correspond to points of the rational
normal curve, while all the points outside the curve have rank
two. Do the same for binary cubics and stratify $\mathbb{P}
S_3=\mathbb{P}^3$.
\end{exercise}

We can produce a useful interpretation of Terracini's Lemma in the
case of Veronese varieties. We consider the Veronese variety
$X=\nu_d(\mathbb{P}^n)\subset\mathbb{P}^N$.

\begin{rem} {\it If $H\subset\mathbb{P}^N$ is a hyperplane, then
$\nu_d^{-1}(H)$ is a degree $d$ hypersurface.} To see this, notice
that $H$ has an equation of the form $a_0z_0+\ldots +a_N z_N$
where $z_i$ are the coordinates of $\mathbb{P}^N$. To determine an
equation for $\nu_d^{-1}(H)$ it is enough to substitute each $z_i$
with the corresponding degree $d$ monomial in the
$x_0,\ldots,x_n$.
\end{rem}

\begin{rem} {\it If $H\subset\mathbb{P}^N$ is a hyperplane and $[L^d]\in H$, then
$\nu_d^{-1}(H)$ is a degree $d$ hypersurface passing through the
point $[L]\in\mathbb{P}^n$.} This is clearly true since
$\nu_d^{-1}([L^d])=[L]$.
\end{rem}

\begin{rem} {\it If $H\subset\mathbb{P}^N$ is a hyperplane such that $T_{[L^d]}(X)\subset H$, then
$\nu_d^{-1}(H)$ is a degree $d$ hypersurface singular at the point
$[L]\in\mathbb{P}^n$.} This can be seen using {\it apolarity} or
by
 direct computation choosing $L^d=x_0^d$.
\end{rem}

We illustrate the last remark in an example.

\begin{example} Consider the Veronese surface $X \subseteq
\mathbb{P}^5$, let $$P=[1:0:0:0:0:0]=[x^2] \in X\text{,}$$ and let
$\mathbb{C}[z_0,z_1,\ldots,z_5]$ be the coordinate ring of
$\mathbb{P}^5$. If $H$ is a hyperplane containing $P$, then $H$ has
equation
\[a_1z_1+a_2z_2+a_3z_3+a_4z_4+a_5z_5 = 0\]
and hence $\nu_2^{-1}(H)$ is the plane conic determined by the equation
\[a_1xy+a_2xz+a_3y^2+a_4yz+a_5z^2=0,\]
which passes through the point $\nu^{-1}(P)=[1:0:0]$. The tangent
space $T_P(X)$ is the projective space associated to the linear span of the forms
\[x^2,xy,xz,\]
and hence it is the linear span of the points
\[[1:0:0:0:0:0],[0:1:0:0:0:0],[0:0:1:0:0:0].\]
Thus, if $H\supset T_P(X)$ then $a_1=a_2=0$ and the corresponding
conic has equation
\[a_3y^2+a_4yz+a_5z^2=0,\]
which is singular at the point $[1:0:0]$.
\end{example}

\begin{exercise}\label{exTWO2} Repeat the argument above to prove the general
statement: if $T_{[L^d]}(\nu_d(\mathbb{P}^n))\subset H$, then
$\nu_d^{-1}(H)$ is a degree $d$ hypersurface singular at the point
$[L]\in\mathbb{P}^n$.
\end{exercise}

We will now elaborate on the connection between double point schemes and
higher secant varieties to Veronese varieties.

\begin{defn} Let $P_1,\ldots,P_s\in\mathbb{P}^n$ be points with defining
ideals $\wp_1,\ldots,\wp_s$ respectively. The scheme defined by
the ideal $\wp_1^2\cap\ldots\cap\wp_s^2$ is called a {\it $2$-fat
point} scheme or a {\it double point} scheme.
\end{defn}

\begin{rem} Let $X=\nu_d(\mathbb{P}^n)\subset\mathbb{P}^N$. There is a bijection
between
\[\lbrace H\subset\mathbb{P}^N \mbox{ a hyperplane }: H \supset \langle T_{P_1}(X),\ldots,T_{P_s}(X)\rangle \rbrace,\]
and
\[\lbrace\mbox{degree $d$ hypersurfaces of }\mathbb{P}^n\mbox{ singular at } P_1,\ldots,P_s\rbrace=(\wp_1^2\cap\ldots\cap\wp_s^2)_d.\]
\end{rem}

Using the double point interpretation of Terracini's Lemma we get
the following criterion to study the dimension of higher secant
varieties to Veronese varieties.

\begin{lem} Let $X=\nu_d(\mathbb{P}^n)\subset\mathbb{P}^N$ and
choose generic points $P_1,\ldots,P_s\in\mathbb{P}^n$ with
defining ideals $\wp_1,\ldots, \wp_s$ respectively. Then
\[\dim\sigma_s(X)=N-\dim (\wp_1^2\cap\ldots\cap\wp_s^2)_d.\]
\end{lem}

\begin{example} We consider, again, the Veronese surface $X$ in
$\mathbb{P}^5$. To determine $\dim\sigma_2(X)$ we choose generic
points $P_1,P_2\in\mathbb{P}^2$ and look for conics singular at
both points, \emph{i.e.} elements in $(\wp_1^2\cap\wp_s^2)_2$. Exactly
one such conic exists (the line through $P_1$ and $P_2$ doubled)
and hence $\sigma_2(X)$ is a hypersurface.
\end{example}

\begin{exercise}\label{exTWO3} Solve the big Waring problem for $n=1$ using the
double points interpretation.
\end{exercise}

We now return to the big Waring problem. Notice that the secant variety
interpretation and a straightforward dimension count yields an
expected value for $g(n,d)$ which is
\[\left\lceil{{d+n \choose n}\over n+1}\right\rceil.\]
This expectation turns out to be true except for a short list of exceptions. 
 A complete solution for the big Waring problem is given by
a celebrated result by Alexander and Hirschowitz, see \cite{AH92}.

\begin{thm}[\cite{AH92}] Let $F$ be a generic degree $d$ form in $n+1$
variables. Then
\[\mathrm{rk}(F)=\left\lceil{{d+n \choose n}\over n+1}\right\rceil,\]
unless
\begin{itemize}
\item $d=2$, any $n$ where $\mathrm{rk}(F)=n+1$.

\item $d=4,n=2$ where $\mathrm{rk}(F)=6$ and not $5$ as expected.

\item $d=4,n=3$ where $\mathrm{rk}(F)=10$ and not $9$ as expected.

\item $d=3,n=4$ where $\mathrm{rk}(F)=8$ and not $7$ as expected.

\item $d=4,n=4$ where $\mathrm{rk}(F)=15$ and not $14$ as
expected.

\end{itemize}
\end{thm}

\begin{rem} A straightforward interpretation of the Alexander and
Hirschowitz result in terms of higher secants is as follows.  The number
$g(n,d)$ is the smallest $s$ such that
$\sigma_s(\nu_d(\mathbb{P}^n))=\mathbb{P}^N$,
unless $n$ and $d$ fall into one of the exceptional cases above.
\end{rem}

\begin{rem} Actually the Alexander and Hirschowitz result gives
more for higher secant varieties of the Veronese varieties, namely
that $\nu_d(\mathbb{P}^n)$ is not defective, for all $s$, except for
the exceptional cases.
\end{rem}

Let's now try to explain some of the defective cases of the
Alexander-Hirschowitz result.

\begin{example} For $n=2,d=4$ we consider $X=\nu_4(\mathbb{P}^2)\subset\mathbb{P}^{14}$.
In particular, we are looking for the smallest $s$ such that
$\sigma_s(X)=\mathbb{P}^{14}$. We expect $s=5$ to work and we want
to check whether this is the case or not. To use the double points
interpretation, we choose $5$ generic points $P_1,\ldots,
P_5\in\mathbb{P}^2$ and we want to determine
\[ \dim (\wp_1^2\cap\ldots\cap\wp_5^2)_4 \text{.}\]
To achieve this, we want to know the dimension of the space of quartic curves that are
singular at each $P_i$. Counting conditions we have $15-5\times
3 = 0$ and expect that
$$ \dim (\wp_1^2\cap\ldots\cap\wp_5^2)_4 = 0 \text{.}$$ In fact, there exists a conic passing
through the points $P_i$ and this conic doubled is a quartic with
the required properties. Thus,
\[\dim (\wp_1^2\cap\ldots\cap\wp_5^2)_4\geq 1,\]
and $\dim\sigma_5(X)\leq 14-1=13$.
\end{example}

\begin{exercise}\label{exTWO4} Show that $\sigma_5(\nu_4(\mathbb{P}^2))$ is a hypersurface, \emph{i.e.} that it has dimension equal to $13$.
\end{exercise}

\begin{exercise}\label{exTWO5} Explain the exceptional cases $d=2$ any $n$.
\end{exercise}

\begin{exercise}\label{exTWO6} Explain the exceptional cases $d=4$ and $n=3,4$.
\end{exercise}

\begin{exercise}\label{exTWO7} Explain the exceptional case $d=3$ and $n=4$.
(Hint: use Castelnuovo's Theorem which asserts that there exists a
(unique) rational normal curve passing through $n+3$ generic
points in $\mathbb{P}^n$.)
\end{exercise}


\section{Lecture Three}\label{lectthree}

In the last lecture we explained the solution to the big Waring problem and showed how to determine the Waring rank
$\mathrm{rk}(F)$ for $F$ a generic form. We now focus on a more
general  question: given any form $F$ what can we say about
$\mathrm{rk}(F)$?

The main tool we will use is {\it Apolarity} and, in order to do
this, we need the following setting. Let
$S=\mathbb{C}[x_0,\ldots,x_n]$ and $T=\mathbb{C}[y_0,\ldots,y_n]$.
We make $T$ act on $S$ via differentiation, \emph{i.e.} we define
\[y_i\circ x_j={\partial \over \partial x_i}x_j,\]
\emph{i.e.} $y_i\circ x_j=1$ if $i=j$ and it is zero otherwise. We then
extend the action to all $T$ so that $\partial\in T$ is seen as a
differential operator on elements of $S$; from now on we will omit
$\circ$.  If $A$ is a subset of a graded ring, we let $A_{d}$ denote the degree $d$ graded piece of $A$.

\begin{defn} Given $F\in S_d$ we define the {\it annihilator}, or
{\it perp ideal}, of $F$ as follows:
\[F^\perp=\{ \partial\in T : \partial F=0\}.\]
\end{defn}

\begin{exercise}\label{exTHREE1} Show that $F^\perp\subset T$ is an ideal and that
it also is Artinian, \emph{i.e.}  $(T/F^\perp)_i$ is zero for $ i > d$.
\end{exercise}

\begin{exercise}\label{exTHREE2} Let $S_i$ and $T_i$ denote the degree $i$ homogenous pieces of $S$ and $T$ respectively.  Show that the map
\[S_i\times T_i \longrightarrow \mathbb{C}\]
\[(F,\partial)\mapsto \partial F\]
is a perfect pairing, \emph{i.e.}
\[(F,\partial_0)\mapsto 0, \forall F\in S_i \Longrightarrow \partial_0=0,\]
and
\[(F_0,\partial)\mapsto 0, \forall \partial\in T_i \Longrightarrow F_0=0.\]
\end{exercise}

\begin{rem}{ Recall that Artinian Gorenstein rings are characterized by the property that they are all of the form $A = T/F^{\perp}$. Moreover, a property of such an $A$ is that it is finite dimensional, and the Hilbert function is symmetric.}
\end{rem}

\begin{rem} {\it Actually even more is true, and $A=T/F^\perp$ is
Artinian and Gorenstein with socle degree $d$.} Using the perfect
pairing $S_i\times T_i \longrightarrow \mathbb{C}$ we see that
$\dim A_d=\dim A_0=1$ and that $A_d$ is the socle of $A$.
\end{rem}

In what follows we will make use of Hilbert functions, thus we
define them here.

\begin{defn} For an ideal $I\subset T$ we define the {\it Hilbert
function} of $T/I$ as
\[HF(T/I,t)=\dim (T/I)_t.\]
\end{defn}

\begin{example} Let $F\in S_d$. We see that
$HF(T/F^\perp,t)=0$ for all $t > d$, in fact all partial differential operators of
degree $t>d$ will annihilate the degree $d$ form $F$ and hence
$(T/F^\perp)_t = 0$, for $t > d$. From the remark above we also see that
$HF(T/F^\perp,d)=1$.
\end{example}

\begin{exercise}\label{exTHREE3} Given $F\in S_d$ show that $HF(T/F^\perp,t)$ is a symmetric
function with respect to ${d+1 \over 2}$ of $t$.
\end{exercise}

An interesting property of the ideal $F^\perp$ is described by
Macaulay's Theorem (see \cite{Mac1927}).

\begin{thm}
If $F\in S_d$, then $T/F^\perp$ is an Artinian Gorenstein ring
with socle degree $d$. Conversely, if $T/I$ is an Artinian
Gorenstein ring with socle degree $d$, then $I=F^\perp$ for some
$F\in S_d$.
\end{thm}

Let's now see how apolarity relates to the Waring rank. Recall
that $s=\mathrm{rk}(F)$ if and only if $F=\sum_1^s L_i^d$ and no
shorter presentation exists.

\begin{example} We now compute the possible Waring ranks for a
binary cubic, \emph{i.e.} for $F\in S_3$ where $S=\mathbb{C}[x_0,x_1]$.
We begin by describing the Hilbert function of $F^\perp$. There
are only two possibilities:
\begin{itemize}
\item[{\bf case 1}]
\[
\begin{array}{llllll}

t               & 0 &   1   &   2   &   3   & 4 \\

\hline

 HF(T/F^\perp,t) & 1 &   1   &   1   &   1   & 0\rightarrow
\end{array}
\]

\item[{\bf case 2}]
\[
\begin{array}{llllll}

t               & 0 &   1   &   2   &   3   & 4 \\

\hline

 HF(T/F^\perp,t) & 1 &   2   &   2   &   1   & 0\rightarrow
\end{array}
\]

\end{itemize}

We want to show that in {\bf case $1$} we have $F=L^3$. From the
Hilbert function we see that $(F^\perp)_1=\langle \partial_1
\rangle$. From the perfect pairing property we see that
 \[ \{L\in S_1 : \partial_1 L=0\}=\langle L_1
\rangle. \]

Thus we can find $L_0\in S_1$ such that $\partial_1 L_0=1$ and
\[S_1=\langle x_0,x_1 \rangle=\langle L_0,L_1 \rangle.\]
We now perform a linear change of variables and we obtain a
polynomial
$$G(L_0,L_1)=aL_0^3+bL_0^2L_1+cL_0L_1^2+dL_1^3$$ such
that
\[G(L_0,L_1)=F(x_0,x_1).\]
As $\partial_1 L_0\neq 0$ and $\partial_1 L_1= 0$ we get
\[0=\partial_1 G=2bL_0L_1+cL_1^2+3dL_1^2,\]
and hence $G=F=aL_0^3$ thus $\mathrm{rk}(F)=1$.

We want now to show that in {\bf case $2$} we have
$\mathrm{rk}(F)= 2$ or $\mathrm{rk}(F)= 3$. We note that
$\mathrm{rk}(F)\neq 1$, otherwise $(F^\perp)_1\neq 0$. As in this
case $(F^\perp)_1=0$, we consider the degree two piece,
$(F^\perp)_2=\langle Q\rangle$. We have to possibilities
\[Q=\partial\partial'\mbox{, where $\partial$ and $\partial'$ are not proportional, or } Q=\partial^2.\]

If $Q=\partial\partial'$, where $\partial$ and $\partial'$ are not proportional, we can construct a basis for $S_1=\langle
L,L'\rangle$ in such a way that
\[\partial L=\partial' L'=1,\]
and
\[\partial' L=\partial L'=0.\]
Then we perform a change of variables and obtain
\[F(x_0,x_1)=G(L_0,L_1)=aL_0^3+bL_0^2L_1+cL_0L_1^2+dL_1^3.\]
We want to show that $F(x_0,x_1)=aL_0^3+dL_1^3$. To do this we
define
\[H(x_0,x_1)=G(L_0,L_1)-aL_0^3-dL_1^3,\]
and show that the degree $3$ polynomial $H$ is the zero
polynomial. To do this, it is enough to show that
$(H^\perp)_3=T_3$. We now compute that
\[\partial^3 H= 6a L - 6 aL=0 \text{ and ,}\]
\[\partial'^3 H= 6d L' - 6 dL'=0 \text{.}\]
We then notice that $\partial^2\partial'=\partial Q\in F^\perp$
and $\partial^2\partial' H=0$; similarly for
$\partial\partial'^2$. Thus $H=0$ and $F(x_0,x_1)=aL_0^3+dL_1^3$.
As $(F^\perp)_1=0$ this means that $\mathrm{rk}(F)=2$.

Finally, if $Q=\partial^2$ we assume by contradiction that
$\mathrm{rk}(F)= 2$, thus $F=N^3+M^3$ for some linear forms $N$
and $M$. There exist linearly
independent differential operators $\partial_N,\partial_M\in S_1$ such that
\[\partial_N N=\partial_M M=1,\]
and
\[\partial_N M=\partial_M N=0.\]
And then $\partial_N\partial_M\in F^\perp$ and this is a
contradiction as $Q$ is the only element in $(F^\perp)_2$ and it
is a square.
\end{example}

\begin{rem} We consider again the case of binary cubic forms. We want to make a connection
between the Waring rank of $F$ and certain ideals contained in
$F^\perp$. If $\mathrm{rk}(F)=1$ then we saw that $F^\perp\supset
(\partial_1)$ and this is the ideal of one point in
$\mathbb{P}^1$. If $\mathrm{rk}(F)=2$ then $F^\perp\supset
(\partial\partial')$ and this the ideal of two distinct points in
$\mathbb{P}^1$; as $(F^\perp)_1=0$ there is no ideal of one point
contained in the annihilator. Finally, if $\mathrm{rk}(F)=3$, then
$F^\perp\supset (\partial^2)$ and there is no ideal of two points,
or one point, contained in the annihilator. However,
$(F^\perp)_3=T_3$ and we can find many ideals of three points.
\end{rem}

There is a connection between $\mathrm{rk}(F)$ and set of points
whose ideal $I$ is such that $I\subset F^\perp$. This connection
is the content of the {\it Apolarity Lemma}, see \cite{IarrobinoKanev_text}.

\begin{lem} Let $F\in S_d$ be a degree $d$ form in $n+1$
variables. Then the following facts are equivalent:
\begin{itemize}
\item $F=L_1^d+\ldots +L_s^d$;

\item $F^\perp\supset I$ such that $I$ is the ideal of a set of
$s$ distinct points in $\mathbb{P}^n$.

\end{itemize}

\end{lem}

\begin{example} We use the Apolarity Lemma to explain the
Alexander-Hirschowitz defective case $n=2$ and $d=4$. Given a
generic $F\in S_4$ we want to show that $\mathrm{rk}(F)=6$ and not
$5$ as expected. To do this we use Hilbert functions. Clearly, if
$I\subset F^\perp$ then $HF(T/I,t)\geq HF(T/F^\perp,t)$ for all
$t$. Thus by computing $HF(T/F^\perp,t)$ we get information on the
Hilbert function of any ideal contained in the annihilator, and in
particular for ideal of sets of points.

\[
\begin{array}{llllllll}

t               & 0 &   1   &   2   &   3   & 4  & 5 & 6\\

\hline

 HF(T/F^\perp,t) & 1 &   3   &   6   &   3   & 1 & 0 & \rightarrow
\end{array}
\]

In particular, $HF(T/F^\perp,2)=6$ means that for no set of $5$
points its defining ideal $I$ could be such that $I\subset
F^\perp$.

\end{example}

\begin{exercise}\label{exTHREE4} Use the Apolarity Lemma to compute
$\mathrm{rk}(x_0x_1^2)$. Then try the binary forms $x_0x_1^d$.
\end{exercise}

\begin{exercise}\label{exTHREE5} Use the Apolarity Lemma to explain
Alexander-Hirschowitz exceptional cases.
\end{exercise}

It is in general very difficult to compute the Waring rank of a
given form and (aside from brute force) no algorithm exists which can compute it in
all cases.  Lim and Hillar show that this problem is an instance of the fact that, as their title states, ``Most tensor problems are NP-Hard,'' \cite{HillarLim09mosttensor}.
However, we know $\mathrm{rk}(F)$ when $F$ is a
quadratic form, and we do have an efficient algorithm when $F$ is
a binary form.

\begin{rem} {\it There is an algorithm, attributed to Sylvester, to
compute $\mathrm{rk}(F)$ for a binary form and it uses the
Apolarity Lemma.} The idea is to notice that
$F^\perp=(\partial_1,\partial_2)$, \emph{i.e.} the annihilator is a
complete intersection ideal, say, with generators in degree
$d_1=\deg\partial_1\leq d_2=\deg\partial_2$. If $\partial_1$ is
square free, then we are done and $\mathrm{rk}(F)=d_1$. If not, as
$\partial_1$ and $\partial_2$ do not have common factors, there is
a square free degree $d_2$ element in $F^\perp$. Hence,
$\mathrm{rk}(F)=d_2$.
\end{rem}

\begin{exercise}\label{exTHREE6} Compute $\mathrm{rk}(F)$ when $F$ is a
quadratic form.
\end{exercise}

\begin{rem} The Waring rank for monomials was determined in 2011 in a
paper of Carlini, Catalisano and Geramita, see \cite{CCG_Waring}, and independently by Buczy{\'n}ska, Buczy{\'n}ski, and Teitler, see \cite{BBT_Waring}. In particular, it was shown that
\[\mathrm{rk}(x_0^{a_0}\ldots x_n^{a_n})={1 \over (a_0+1)}\Pi_{i=0}^n (a_i+1),\]
where $1\leq a_0\leq a_1\leq\ldots\leq a_n$.
\end{rem}

We conclude this lecture by studying the Waring rank of degree $d$ forms of the
kind $L_1^d+\ldots + L_s^d$. Clearly, $\mathrm{rk}(L_1^d)=1$ and
$\mathrm{rk}(L_1^d+L_2^d)=2$, if $L_1$ and $L_2$ are linearly
independent. If the linear forms $L_i$ are not linearly
independent, then the situation is more interesting.

\begin{example} Consider the binary cubic form
$F=ax_0^3+bx_1^3+(x_0+x_1)^3$. We want to know $\mathrm{rk}(F)$.
For a generic choice of $a$ and $b$, we have $\mathrm{rk}(F)=2$,
but for special values of $a$ and $b$ $\mathrm{rk}(F)=3$. The idea
is that the rank three element of $\mathbb{P}S_3$ lie on the
tangent developable of the twisted cubic curve, which is an
irreducible surface. Hence, the general element of the plane
\[\langle [x_0^3],[x_1^3],[(x_0+x_1)^3]\rangle\]
has rank two, but there are rank three elements.
\end{example}

\begin{exercise}\label{exTHREE7} Prove that $\mathrm{rk}(L^d+M^d+N^d)=3$ whenever
$L,M$ and $N$ are linearly independent linear forms.
\end{exercise}


\section{Lecture Four}\label{lectfour}

In the last lecture we introduced the Apolarity Lemma and used
it to study the Waring rank of a given specific form. In this
lecture we will go back to the study of higher secant varieties of
varieties that are not Veronese varieties.

The study of higher secant varieties of Veronese varieties is
connected to Waring's problems, and hence with sum of powers
decompositions of forms. We now want to consider tensors in general
and not only homogenous polynomials, which correspond to symmetric tensors.

Consider $\CC$-vector spaces $V_1,\ldots,V_t$ and the tensor product
\[V=V_1\otimes\ldots\otimes V_t.\]

\begin{defn} A tensor $v_1\otimes\ldots\otimes v_t\in V$ is called {\it elementary}, or \emph{indecomposable} or \emph{rank-one} tensor.
\end{defn}
Elementary tensors are the building blocks of $V$. More specifically, there is a basis of $V$ consisting of elementary tensors, so any tensor
$T\in V$ can be written as a linear combination of elementary tensors;
in this sense elementary tensors are analogous to monomials.

A natural question is: given a tensor $T$ what
is the minimum $s$ such that $T=\sum_1^s T_i$ where each $T_i$ is
an elementary tensor? The value $s$ is called the {\it tensor rank} of $T$ and is the analogue of Waring rank for forms. Of course we could state tensor versions of the Waring's problems and try to solve them as well.

In order to study these problems geometrically, we need to
introduce a new family of varieties.

\begin{defn} Given vector spaces $V_1,\ldots, V_t$ the {\it Segre
map} is the map
\[\mathbb{P}V_1\times\ldots\times\mathbb{P}V_t\longrightarrow\mathbb{P}(V_1\otimes\ldots\otimes V_t)\]
\[([v_1],\ldots,[v_t])\mapsto [v_1\otimes\ldots\otimes v_t] \]
and the image variety  is called a {\it Segre variety} or {\it
the Segre product} of $\mathbb{P}V_1,\ldots,\mathbb{P}V_t$.
\end{defn}

If the vector spaces are such that $\dim V_i=n_i+1$, then we will
often denote by
\[X=\mathbb{P}^{n_1}\times\ldots\times\mathbb{P}^{n_t}\]
the image of the Segre map.   In particular,
$X\subset\mathbb{P}(V_1\otimes\ldots\otimes V_t)=\mathbb{P}^N$
where $N+1=\Pi (n_i+1)$. Note that $\dim X=n_1+\ldots+n_t$.

By choosing bases of the vector spaces $V_i$ we can write the
Segre map in coordinates
\[\mathbb{P}^{n_1}\times\mathbb{P}^{n_2}\times\ldots\times\mathbb{P}^{n_t}\longrightarrow\mathbb{P}^N \]
\[
\begin{array}{l}
[a_{0,1}:\ldots:a_{n_1,1}]\times[a_{0,2}:\ldots:a_{n_2,2}]\times\ldots\times[a_{0,t}:\ldots:a_{n_t,t}]\mapsto

\\
\phantom{xxxxxxx} [a_{0,1}a_{0,2}\ldots
a_{0,t}:a_{0,1}a_{0,2}\ldots
a_{1,t}:\ldots:a_{n_1,1}a_{n_2,2}\ldots a_{n_t,t}]
\end{array}.
\]

\begin{example} Consider
$X=\mathbb{P}^{1}\times\mathbb{P}^{1}$, then
$X\subset\mathbb{P}^3$ is a surface. The Segre map is
\[\mathbb{P}^{1}\times\mathbb{P}^{1}\longrightarrow\mathbb{P}^3\]
\[[a_0:a_1]\times[b_0:b_1]\mapsto[a_0b_0:a_0b_1:a_1b_0:a_1b_1].\]
If $z_0,z_1,z_2$, and $z_3$ are the coordinates of
$\mathbb{P}^3$, then it is easy to check that $X$ has equation
\[z_0z_3-z_1z_2=0.\]
Thus $X$ is a smooth quadric in $\mathbb{P}^3$.
\end{example}

\begin{example} Consider again
$X=\mathbb{P}^{1}\times\mathbb{P}^{1}\subset\mathbb{P}^3$. We identify $\mathbb{P}^3$ with the projectivization of the vector space of $2\times 2$ matrices. Using this identification we can
write the Segre map as
\[\mathbb{P}^{1}\times\mathbb{P}^{1}\longrightarrow\mathbb{P}^3\]
\[[a_0:a_1]\times[b_0:b_1]\mapsto
\left(\begin{array}{cc}
a_0b_0 & a_0b_1 \\
a_1b_0 & a_1b_1
\end{array}\right).\] Thus $X$
represents the set of $2\times 2$ matrices of rank at most one and
the ideal of $X$ is generated by the vanishing of the determinant of the
generic matrix $\left(\begin{array}{cc}z_0 & z_1 \\ z_2 &
z_3\end{array}\right).$
\end{example}

\begin{exercise}\label{exFOUR1} Workout a matrix representation for the Segre
varieties with two factors
$\mathbb{P}^{n_1}\times\mathbb{P}^{n_2}$.
\end{exercise}

Before entering into the study of the higher secant varieties of Segre
varieties, we provide some motivation coming from Algebraic
Complexity Theory.

\begin{example} The multiplication of two $2\times 2$ matrices can
be seen as bilinear map
\[
T\colon \CC^{4}\times \CC^{4}\to \CC^{4}
\]
or, equivalently, as a tensor
\[T\in{\mathbb{C}^4}^*\otimes{\mathbb{C}^4}^*\otimes\mathbb{C}^4.\]
It is interesting to try to understand how many multiplications over the ground field are required to compute the map $T$.

If we think of $T$ as a tensor, then we can write it as a linear combination of elementary tensors and each elementary tensor represents a multiplication.
The naive algorithm for matrix multiplication, which in general uses $n^{3}$ scalar multiplications to compute the product of two $n\times n$ matrices,  implies
that $T$ can be written as $T=\sum_1^8\alpha_i\otimes\beta_i\otimes
c_i$. However, Strassen in \cite{Strassen83_rank} proved that
\[[T]\in\sigma_7(\mathbb{P}^3\times\mathbb{P}^3\times\mathbb{P}^3),\]
and, even more, that $T$ is the sum of $7$
elementary tensors
\[T=\sum_1^7\alpha_i\otimes\beta_i\otimes
c_i.\]

Strassen's algorithm actually holds for multiplying matrices over any algebra.  Thus by viewing a given $n\times n$ matrix in one with size a power of $2$, one can use Strassen's algorithm iteratively. So $2^{m}\times 2^{m}$ matrices can be multiplied using $7^{m}$ multiplications.
In particular, this method lowers the upper bound for the complexity of matrix multiplication from $n^3$ to $n^{\log_{2}7}\simeq n^{2.81}$.

After Strassen's result, it was shown that the rank of $T$ is not smaller than $7$.
On the other hand, much later, Landsberg proved that the border rank of $T$ is $7$, that is to say that $T \not\in \sigma_6(\mathbb{P}^3\times\mathbb{P}^3\times\mathbb{P}^3)$, \cite{Lan_2by2_matrix}.
\end{example}

The question of the complexity of matrix multiplication has recently been called one of the most important questions in Numerical Analysis \cite{SIAM_news}.  The reason for this is that the complexity of matrix multiplication also determines the complexity of matrix inversion.  Matrix inversion is one of the main tools for solving a square system of linear ODE's.

Williams, in 2012, improved the Coppersmith-Winograd algorithm to obtain the current best upper bound for the complexity of matrix multiplication \cite{Williams12}, but it would lead us to far afield to discuss this here.  On the other hand, the current best lower bounds come from the algebra and geometry of secant varieties.  These bounds arise by showing non-membership of the matrix multiplication tensor on certain secant varieties. To do this, one looks for non-trivial equations that vanish on certain secant varieties, but do not vanish on the matrix multiplication tensor. Indeed, the best results in this direction make use of representation theoretic descriptions of the ideals of secant varieties. For more, see \cite{Landsberg_NewLower,LanOtt_NewLower}.

\subsection{Dimension of secant varieties of Segre varieties}
In the $2\times 2$ matrix multiplication example, we should have pointed out the fact that $\sigma_{7}(\PP^{3}\times \PP^{3}\times \PP^{3})$ actually fills the ambient space, so almost all tensors in $\PP(\CC^{4}\otimes \CC^{4}\otimes \CC^{4})$ have rank 7.  For this and many other reasons we would like to know the dimensions of secant varieties of Segre varieties.

Like in the polynomial case, there is an expected dimension, which is obtained by the naive dimension count.  When $X$ is the Segre product $\mathbb{P}^{n_1}\times\ldots\times\mathbb{P}^{n_t}$, the expected dimension of $\sigma_{s}(X)$ is
\[
\min\left\{s\sum_{i=1}^{t}(n_{i}+1) +s -1  ,\prod_{i=1}^{t}(n_{i}+1)-1\right\}
.\]

We would like to have an analogue of the Alexander-Hirschowitz theorem for the Segre case, however, this is a very difficult problem.  See \cite{AOP_Segre} for more details.  There are some partial results, however. For example Catalisano, Geramita and Gimigliano in \cite{CGG_P1s} show that $\sigma_{s}(\PP^{1}\times \dots \times \PP^{1})$ always has the expected dimension, except for the case of $4$ factors.

Again, the first tool one uses to study the dimensions of secant varieties is the computation of tangent spaces together with Terracini's lemma.

\begin{exercise}\label{exFOUR2}
Let $X=\PP V_{1}\times \dots \times \PP V_{t}$ and let $[v] =[v_{1}\otimes\dots\otimes v_{t}]$ be a point of $X$.
Show that the cone over the tangent space to $X$ at $v$ is the span of the following vector spaces:
\[\begin{matrix}
V_{1}\otimes v_{2}\otimes v_{3}\otimes \dots \otimes v_{t},\\
v_{1}\otimes V_{2}\otimes v_{3}\otimes \dots \otimes v_{t},\\
\vdots \\
v_{1}\otimes v_{2}\otimes \dots \otimes v_{t-1}\otimes V_{t}.
\end{matrix}\]
\end{exercise}

\begin{exercise}\label{exFOUR3}
Show that $\sigma_{2}(\PP^{1}\times \PP^{1}\times \PP^{1}) = \PP^{7}$.
\end{exercise}

\begin{exercise}\label{exFOUR4}
Use the above description of the tangent space of the Segre product and Terracini's lemma to show that $\sigma_{3}(\PP^{1}\times \PP^{1}\times \PP^{1}\times \PP^{1})$ is a hypersurface in $\PP^{15}$ and not the entire ambient space as expected. This shows that the four-factor Segre product of $\PP^{1}$s is defective.
\end{exercise}
There are two main approaches to the study of the dimensions of
secant varieties of Segre products: \cite{CGG3_SegreVeronese} and
\cite{AOP_Segre}. In  \cite{CGG3_SegreVeronese} the authors
introduce and use what they call the affine-projective method. In this
way, the study of the dimension of higher secant varieties of
Segre, and Segre-Veronese, varieties reduces to the study of the
postulation of non-reduced schemes supported on linear spaces. In
\cite{AOP_Segre},  the authors show that the ``divide and
conquer'' method of Alexander and Hirschowitz can be used to set
up a multi-step induction proof for the non-defectivity of Segre
products.  They are able to obtain partial results on
non-defectivity by then checking many initial cases, often using
the computer.  On the other hand, for the remaining cases there
are many more difficult computations to do in order to get the
full result.

\subsection{Flattenings}
Often, the first tool used to understand properties of tensors is to reduce to Linear Algebra (when possible). For this, the notion of flattenings is essential.
Consider for the moment the three-factor case. We may view the vector space $\CC^{a}\otimes \CC^{b}\otimes \CC^{c}$ as a space of matrices in three essentially different ways as the following spaces of linear maps:
\[
\begin{matrix}
(\CC^{a})^{*} \to \CC^{b} \otimes \CC^{c} \\
(\CC^{b})^{*} \to \CC^{a} \otimes \CC^{c}\\
(\CC^{c})^{*} \to \CC^{a} \otimes \CC^{b}
\end{matrix}
.\]
(\emph{A priori} there are many more choices of flattenings, however, in the three-factor case the others are obtained by transposing the above maps.)
When there are more than 3 factors the situation is similar, with many more flattenings to consider.

For $V_{1}\otimes \dots\otimes V_{t}$, a \emph{$p$-flattening} is the interpretation as a space of matrices with $p$ factors on the left:
\[
(V_{i_{1}}\otimes \dots \otimes V_{i_{p}})^{*} \to V_{j_{1}}\otimes \dots \otimes V_{j_{t-p}}
.\]
For a given tensor $T \in V_{1}\otimes \dots\otimes V_{t}$ we call a $p$-flattening of $T$ a realization of $T$ in one of the above flattenings.  This naturally gives rise to the notion of \emph{multi-linear rank}, which is the vector the ranks of the 1-flattenings of $T$, see \cite{MR2776439}.
\begin{exercise}\label{exFOUR5}
Show that $T$ has rank 1 if and only if its multilinear rank is $(1,\dots,1)$.
\end{exercise}

Recall that a linear mapping $T\colon (\CC^{a})^{*}\to \CC^{b}$ has rank $r$ if the image of the map has dimension $r$ and the kernel has dimension $a-r$.
Moreover, since the rank of the transpose is also $r$, after re-choosing bases in $\CC^{a}$ and $\CC^{b}$ one can find $r$-dimensional subspaces in $\CC^{a}$ and $\CC^{b}$ so that $T\in (\CC^{r})^{*}\otimes \CC^{r}$.
This notion generalizes to tensors of higher order.
In particular, it is well known that  $T\in V_{1}\otimes \dots\otimes V_{t} $ has multilinear rank $\leq (r_{1}  ,\dots,r_{t})$ if an only if there exist subspaces $\CC^{r_{i}}\subset V_{i}$ such that $T\in \CC^{r_{1}}\otimes \dots \otimes \CC^{r_{t}}$.
The Zariski closure of all tensors of multilinear rank $(r_{1},\dots, r_{t})$ is known as the subspace variety, denoted $\Sub_{r_{1},\dots,r_{t}}$.  See \cite{LanWey_Seg} for more details.

The connection between subspace varieties and secant varieties is the content of the following exercise.
\begin{exercise}\label{exFOUR6}\label{sub}
Let $X=\PP V_{1}\times \dots \times \PP V_{t}$.
Show that if $r \leq r_{i}$ for $1\leq i \leq t$, then
\[\sigma_{r}(X)\subset \Sub_{r_{1},\dots,r_{t}}.\]
\end{exercise}
Notice that for the 2-factor case, $\sigma_{r}(\PP^{a-1}\times \PP^{b-1}) = \Sub_{r,r}$.

Aside from the case of binary tensors (tensor products of $\CC^{2}$s), another case that is well understood is the case of very unbalanced tensors. (For the more refined notion of ``unbalanced'' see \cite[\S~4]{AOP_Segre}.)
  Again consider $\CC^{a}\otimes \CC^{b}\otimes \CC^{c}$  and suppose that $a\geq bc$.  Then for all $r\leq bc$ we have  \[\sigma_{r}(\PP^{a-1}\times \PP^{b-1}\times \PP^{c-1}) = \sigma_{r}(\PP^{a-1}\times \PP^{bc-1}) = \Sub_{r,r}.\]
Therefore we can always reduce the case of very unbalanced 3-fold tensors to the case of matrices.  More generally a tensor in $V_{1}\otimes\dots\otimes  V_{t}$ is called very unbalanced if $(n_{i}+1)\geq \prod_{j\neq i}(n_{j}+1)$.  In the very unbalanced case we can reduce to the case of matrices and use results and techniques from linear algebra.  

\subsection{Equations of Secant varieties of Segre products}
Now we turn to the question of defining equations. Recall that a matrix has rank $\leq r$ if and only if all of its $(r+1)\times (r+1)$ minors vanish.  Similarly, a tensor has multilinear rank $\leq (r_{1},\dots, r_{t})$ if all of its $(r_{i}+1)\times (r_{i}+1)$ minors vanish.
It is easy to see that a tensor has rank 1 if and only if all of the $2\times 2$ minors of flattenings vanish.  Since the Segre product is closed, these equations actually define its ideal. For $\sigma_{2}(X)$ it was conjectured by Garcia, Stillman and Sturmfels (GSS) that the $3\times 3$ minors suffice to define the ideal.  The first partial results were by GSS themselves as well as by Landsberg and Manivel \cite{LanMan04_Seg, LanMan08_Strassen}, Landsberg and Weyman \cite{LanWey_Seg} and by Allman and Rhodes \cite{AllmanRhodes03,AllmanRhodes08}, while the full conjecture was resolved by Raicu \cite{RaicuGSS}.

While the result of Exercise~\ref{sub} implies that minors of flattenings give some equations of secant varieties, often they provide no information at all.  For example, consider the case of $3\times 3\times 3$ tensors.  One expects that $\sigma_{4}(\PP^{2}\times \PP^{2}\times \PP^{2})$ fills the entire ambient $\PP^{26}$, however this is not the case.  On the other hand, one cannot detect this from flattenings since all of the flattenings are $3\times 9$ matrices, which have maximum rank $3$, and there are no $5\times 5$ minors to consider.

Strassen noticed that a certain equation actually vanishes on $\sigma_{4}(\PP^{2}\times \PP^{2}\times \PP^{2})$, and moreover, he shows that it is a hypersurface.  Strassen's equation was studied in more generality by \cite{LanMan08_Strassen}, put into a broader context by Ottaviani \cite{Ottaviani09_Waring} generalized by Landsberg and Ottaviani \cite{LanOtt11_Equations}.  Without explaining the full generality of the construction, we can describe Ottaviani's version of Strassen's equation as follows.

Suppose $T\in V_{1}\otimes V_{2}\otimes V_{3}$ with $V_{i}\cong\CC^{3}$ and consider the flattening
\[
(V_{1})^{*}\to V_{2}\otimes V_{3}
.\]
Choose a basis $\{v_{1},v_{2},v_{3}\}$ for $V_{1}$ and write $T$ as a linear combination of matrices $T = v_{1}\otimes T^{1} + v_{2}\otimes T^{2} + v_{3}\otimes T^{3}$. The $3\times 3$ matrices  $T^{i}$ are called the slices of $T$ in the $V_{1}$-direction with respect to the chosen basis.  Now consider the following matrix:
\[\varphi_{T} =
\begin{pmatrix}
0 & T^{1} & -T^{2} \\
-T^{1} & 0 &  T^{3} \\
T^{2} & -T^{3} &0
\end{pmatrix}
,\]
where all of the blocks are $3\times 3$.
\begin{exercise}\label{exFOUR7}
\begin{enumerate}
\item Show that if $T$ has rank $1$ then $\varphi_{T}$ has rank 2.
\item Show that $\varphi$ is additive in its argument, \emph{i.e.} show that $\varphi_{T+T'} = \varphi_{T}+ \varphi_{T'}.$
\end{enumerate}
\end{exercise}
The previous exercise together with the subadditivity of matrix rank implies that if $T$ has tensor rank $r$ then $\varphi_{T}$ has matrix rank $\leq 2r$.  In particular, if $T$ has tensor rank $4$, the determinant of $\varphi_{T}$ must vanish. Indeed $\det(\varphi_{T})$ is Strassen's equation, and it is the equation of the degree 9 hypersurface $\sigma_{4}(\PP^{2}\times \PP^{2}\times \PP^{2}).$
\begin{rem}
This presentation of Strassen's equation $\det(\varphi_{T})$ is very compact yet its expansion in monomials is very large, having 9216 terms.
\end{rem}

This basic idea of taking a tensor and constructing a large matrix whose rank depends on the rank of $T$ is at the heart of almost all known equations of secant varieties of Segre products --- see \cite{LanOtt11_Equations}.  One exception is that of the degree 6 equations in the ideal of $\sigma_{4}(\PP^{2}\times \PP^{2}\times \PP^{3})$.  The only known construction of these equations comes from representation theoretic considerations.  For more details see \cite{OedingBates, LanMan04_Seg}.

Despite this nice picture, we actually know surprisingly little about the defining equations of secant varieties of Segre products in general, and this is an ongoing area of current research.

\section{Solution of the exercises}\label{soltns}

In what follows if $S$ is a subset of $\PP^n$, then $\overline{S}$ denotes the smallest closed subset of $\PP^n$ containing $S$ with the reduced subscheme structure.  On the other hand, $\langle S \rangle$ denotes the smallest linear subspace of $\PP^n$ containing $S$.


\subsection*{Exercise \ref{exONE1}}  Show that the union of chords (secant lines) to a plane conic is closed.  However, the union of the chords of the twisted cubic curve in $\PP^3$ is not.
\begin{proof}[Solution]
Let $X\subseteq \PP^2$ be a plane conic.  It suffices to show that $\cup_{p,q\in X} \langle p,q \rangle = \PP^2$.  If $y\in \PP^2\backslash X$ then there exists a line containing $y$ and intersecting $X$ in two distinct points $p$ and $q$.  Thus, $y\in \langle p,q\rangle$ so $\cup_{p,q\in X} \langle p,q \rangle = \PP^2$.

Let $X\subseteq \PP^3$ be the twisted cubic.  Exercise \ref{exONE2++} implies that $\sigma_2(X)=\PP^3$.  On the other hand, direct calculation shows that the point $[0:1:0:0]\in \PP^3$ lies on no secant line of $X$.  Hence, $\cup_{p,q\in X} \langle p,q \rangle $ is not equal to its closure and hence is not closed.
\end{proof}

\subsection*{Exercise \ref{exONE2}}  For $X\subset\mathbb{P}^N$, show that, if
$\sigma_i(X)=\sigma_{i+1}(X)\neq\mathbb{P}^N$, then $\sigma_i(X)$
is a linear space and hence $\sigma_j(X)=\sigma_i(X)$ for all
$j\geq i$.

\begin{proof}[Solution]
It suffices to prove that, for $k\geq 1$, if $\sigma_k(X)=\sigma_{k+1}(X)$ then $\sigma_{k'}(X)=\langle X \rangle$ for $k'\geq k$.

Note \begin{equation}\label{a} \sigma_{k+1}(X)=\overline{\cup_{p_i,q\in X} \langle \langle p_1,\dots, p_k \rangle, q \rangle} \end{equation}
while $\sigma_k(X) = \overline{\cup_{p_i\in X} \langle p_1,\dots, p_k \rangle}$.  Since the singular locus of $\sigma_k(X)$ is a proper closed subset of $\sigma_k(X)$ there exists a non-singular point $z\in \sigma_k(X)$ such that $z\in \cup_{p_i \in X} \langle p_1,\dots, p_k \rangle$.

Now for all $y\in X$ the line $\langle y,z\rangle$ is contained in $\overline{\cup_{x\in X} \langle x, z \rangle}$ and passes through $z$.  Hence
$\overline{\cup_{x\in X} \langle x,z\rangle} \subseteq T_z \overline{\cup_{x\in X} \langle x,z \rangle}$.   Using \eqref{a}, we deduce
$$X\subseteq  \overline{\cup_{x\in X} \langle x,z\rangle} \subseteq T_z \overline{\cup_{x\in X} \langle x,z \rangle}\subseteq T_z \sigma_{k+1}(X).$$  Hence $\langle X\rangle \subseteq T_z \sigma_{k+1}(X)$.  In addition, since $\sigma_{k+1}(X)=\sigma_k(X)$ we deduce $$\sigma_k(X)\subseteq \langle X \rangle \subseteq T_z \sigma_k(X).$$ Since $z$ is non-singular, $\dim T_z\sigma_k(X)=\dim \sigma_k(X)$.  Finally, since $T_z \sigma_k(X)$ is irreducible and $\sigma_k(X)$ is reduced we conclude $\sigma_k(X)=\langle X \rangle = T_z\sigma_k(X)$.  If $k'\geq k$ then $\langle X \rangle \subseteq \sigma_k(X)\subseteq \sigma_{k'}(X)$.  Since $\sigma_{k'}(X) \subseteq \langle X \rangle$ we deduce $\sigma_{k'}(X)=\sigma_k(X)$.
\end{proof}

\subsection*{Exercise \ref{exONE2+}} If $X\subseteq \PP^N$ is non-degenerate then there exists an $r \geq 1$ with the property that   \[X=\sigma_1(X)\subsetneq\sigma_2(X)\subsetneq \ldots\subsetneq\sigma_r(X)=\mathbb{P}^N.\]
In particular, all inclusions are strict and there is a higher
secant variety that coincides with the ambient space.

\begin{proof}[Solution]
In the notation of Exercise \ref{exONE2}, set $k_0 :=\min \{k\mid \sigma_k(X)=\langle X\rangle\}$.  It suffices to prove that there exists the following chain of strict inclusions
$$X= \sigma_1(X) \subsetneq \sigma_2(X) \subsetneq \dots \subsetneq \sigma_{k_0}(X) = \langle X \rangle.$$
If $\sigma_k(X)=\sigma_{k+1}(X)$ then, by Exercise \ref{exONE2}, $\sigma_k(X)=\langle X \rangle$ and hence $k\geq k_0$.
\end{proof}

\subsection*{Exercise \ref{exONE2++}}  Let $X\subseteq \PP^n$ be a curve.  Prove that $\sigma_2(X)$ has dimension $3$ unless $X$ is contained in a plane.
\begin{proof}[Solution] Let $\mathfrak{X}:=\{ (\langle p, q \rangle , y ) : y \in \langle p, q \rangle \}\subseteq \mathbb{G}(1,3)\times \PP^n$ be the incident correspondence corresponding to the (closure) of the secant line map $(X\times X) - \Delta_X \rightarrow \mathbb{G}(1,3)\times \PP^n$.  Then $\mathfrak{X}$ is an irreducible closed subset of $\mathbb{G}(1,3)\times \PP^n$ and $\sigma_2(X)=p_2(\mathfrak{X})$.  Hence, $\dim \sigma_2(X)\leq 3$.  Suppose now that $\dim \sigma_2(X)=2$.  Let us prove that $X$ is contained in a plane.  First note that for a fixed $p \in X$, $\overline{\cup_{q\in X} \langle p, q\rangle}$ is reduced, irreducible, has dimension $\dim X +1$ and is contained in $\sigma_2(X)$.  Hence, if $\sigma_2(X)$ has dimension $\dim X+1$, $\sigma_2(X)=\overline{\cup_{q\in X} \langle p,q \rangle}$.  Hence, there exists a non-singular $x\in \sigma_2(X)$ such that $x\in \cup_{q \in X} \langle p,q\rangle$ and such that
$ X\subseteq \overline{\cup_{q\in X} \langle x,q\rangle} =\overline{ \cup_{q\in X} \langle p,q\rangle}\subseteq T_x \overline{\cup_{q\in X} \langle x, q \rangle}= T_x \sigma_2(X)$.  Hence, $\langle X \rangle \subseteq T_x \sigma_2(X)$.  Since $\sigma_2(X)$ has dimension $2$ then $T_x \sigma_2(X)$ is a plane.  Finally, note that if $\sigma_2(X)$ has dimension $1$ then $\sigma_2(X)=\sigma_1(X)$.  Hence $\sigma_2(X)=\langle X \rangle$ so $\sigma_2(X)$ and hence $X$ is a line.
\end{proof}

\subsection*{Exercise \ref{exONE2+++}}
Let $M$ be an $n\times n$ symmetric matrix of rank $r$.  Prove that $M$ is a sum of $r$ symmetric matrices of rank $1$.

\begin{proof}[Solution] By performing elementary row and column operations it is possible to find an invertible $n\times n$ matrix $P$ and complex numbers $\lambda_1,\dots, \lambda_n$ such that $$M= P \operatorname{diag}(\lambda_1,\dots, \lambda_n) P^T.$$  Moreover, if $M$ has rank $r$ then there exists $\{i_1,\dots, i_r\}\subseteq \{1,\dots, n\}$ such that $\lambda_i \not = 0$ for $i\in \{i_1,\dots, i_r\}$ and $\lambda_i = 0$ for $i\not \in \{i_1,\dots, i_r\}$.  For $k=1,\dots, n$, set $m^k_{ij}:=\lambda_k p_{ik}p_{jk}$ and let $M_k$ be the matrix with $i,j$ entry $m^k_{ij}$.  Since $m_{ij}^k=m_{ji}^k$,  $M_k$ is symmetric.  Moreover,  $M=\sum_{i=1}^n M_k$.  It remains to show that $M_k$ has rank at most $1$.  Let $P_k$ be the matrix with $i,j$ entry $p_{ik}p_{jk}$.  Since $M_k$ is a scalar multiple of $P_k$ it suffices to show that $P_k$ has rank at most $1$.  For this, we show that all $2\times 2$ minors of $P_k$ are zero.  Indeed an arbitrary $2\times 2$ minor of $P_k$ is the determinant of the $2\times 2$ matrix formed by omitting all rows except row $\alpha, \beta$ and columns $\gamma, \delta$, that is, the determinant $$\left| \begin{matrix} p_{\alpha k}p_{\gamma k} & p_{\alpha k} p_{\delta k} \\ p_{\beta k} p_{\gamma k} & p_{\beta k} p_{\delta k}\end{matrix}\right| $$ which is zero.
\end{proof}

\subsection*{Exercise \ref{exONE3}} Let $X=\nu_2(\PP^2)$.  Let $H$ be the locus of $3\times 3$ symmetric matrices of rank at most $2$.  Prove that $H=\sigma_2(X)$.
\begin{proof}[Solution]  Using Exercise \ref{exONE2+++} we deduce that $X$ is the locus of $3\times 3$ symmetric matrices of rank at most $1$.  On the other hand we know that $X\subseteq \sigma_2(X)\subseteq H$.  Moreover, if $M$ is a symmetric matrix of rank $2$ then, by Exercise \ref{exONE2+++}, $M$ is a sum of two symmetric matrices of rank $1$.  Hence $M\in \sigma_2(X)$.
\end{proof}

\subsection*{Exercise \ref{exONE4}}  Let $H$ be the locus of $4\times 4$ symmetric matrices of rank $2$.  Let $X=\nu_2(\PP^3)$.  Prove that $H=\sigma_2(X)$.  Is $X$ 2-defective?
\begin{proof}[Solution] Using Exercise \ref{exONE2+++} we deduce that $X$ is the locus of $4\times 4$ symmetric matrices of rank at most $1$.  On the other hand we know that $X\subseteq \sigma_2(X)\subseteq H$.  Moreover, if $M$ is a symmetric matrix of rank $2$ then, by Exercise \ref{exONE2+++}, $M$ is a sum of two symmetric matrices of rank $1$.  Hence $M\in \sigma_2(X)$.  To see that $X$ is $2$-defective note that $\operatorname{expdim}(\sigma_2(X))=7$ while the locus of $4\times 4$ symmetric matrices of rank at most $2$ have dimension $6$.
\end{proof}

\subsection*{Exercise \ref{exONE5}}  Prove that if $H\supset T_{P}(X)$, then the
polynomial defining $\nu_{3}^{-1}(H)$ has a double root.
\begin{proof}[Solution]
This is a special case of the following more general claim which is also Exercise \ref{exTWO2}.

\begin{bf} Claim: \end{bf} Let $L\in S_1$ be a liner form.  If $\Lambda$ is a hyperplane in $\PP^N$ containing $T_{[L^d]}(\nu_d(\PP^n))$ then $\nu_d^{-1}(\Lambda)$ is a degree $d$ hypersurface singular at the point $[L]\in \PP^n$.
\begin{proof}
We can find a basis for $S_1$ of the form $\{L,L_1,\dots, L_n\}$ where $L_i$ are linear forms.  With respect to this basis, in coordinates, and using the notation Lecture 1, the image of $[L]$ in $\PP^N$ is the point $p=[1:0:\dots:0]$.  For $i=0,\dots, n$ let $p_i$ be the point of $\PP^N$ with homogeneous coordinates $[Z_0:\dots:Z_N]$ given by $Z_j=\delta_{ij}$.  Translating the intrinsic description of $T_{[L^d]}\nu_d(\PP^n)$ as the linear space $\langle L^{d-1}M \mid M\in S_1\rangle$   into our coordinates for $\PP^n$ and $\PP^N$ with respect to the basis $\{L,L_1,\dots, L_n\}$ for $S_1$ we conclude that $T_P(\nu(\PP^n))=\langle p_0,\dots, p_n\rangle$.  Let $F_{\Lambda}$ be a linear form defining $\Lambda$.  Then $F_{\Lambda}=a_0Z_0+\dots a_N Z_N$.  If $T_p(\nu_d(\PP^n))\subseteq \Lambda$ then $F_{\Lambda}(p_i)=0$, for $i=0,\dots, n$.  Hence $0=a_0=\dots = a_n$.  Let $f$ be the pullback of $F_{\Lambda}$ for $\PP^n$.  Then $f=a_{n+1}x_0^{d-2}x_1+\dots + a_Nx_n^d$.  Clearly the partial derivatives of $f$ vanish at $q=[1:0:\dots :0]$ which is identified with $[L]$.  Hence the zero locus of $f$, a degree $d$ hypersurface, is singular at $q$.  On the other hand, the zero locus of $f$ is $\nu_d^{-1}(\Lambda)$.
\end{proof}

 \end{proof}

\subsection*{Exercise \ref{exONE6}}
Consider the rational normal curve in
$\mathbb{P}^3$, \emph{i.e.} the twisted cubic curve $X=\nu_3(\mathbb{P}(
S_1))\subset\mathbb{P} (S_3)$. We know that $\sigma_2(X)$ fills up
all the space. Can we write any binary cubic as the sum of two
cubes of linear forms? Try $x_0x_1^2$.
\begin{proof}[Solution]  Direct calculation shows that we cannot write $x_0x_1^2$ as a linear combination of cubes of two linear forms.
\end{proof}

\subsection*{Exercise \ref{exONE7}}  Let $X=\nu_d(\PP^n)$.  Recall that $[F]\in X$ if and only if $[F]=L^d$ for some linear form $L$ on $\PP^n$.  Use this description and standard differential geometry to compute $T_{[L^d]}(X)$.
\begin{proof}[Solution]   Let $L$ be a linear form.  Consider an affine curve passing through $L$.  It will have the form $(L+tM)$ where $M$ is allowed to be any linear form.  Considering the image in $\PP^N$ we have
$$(L+tM)^d= L^d+ \binom{d}{d-1}L^{d-1}tM + \text{terms containing higher powers of $t$}.$$  Taking the derivative with respect to $t$ and setting $t=0$ we deduce that
$$T_{[F]}= \langle [L^{d-1}M]\mid M \in S_1\rangle\text{.}$$
\end{proof}


\subsection*{Exercise \ref{exTWO1}}  For binary forms, we can stratify $\PP S_2$ using the Waring rank: rank one elements correspond to points of the rational normal curve, while all the points outside the curve have rank two.  Do the same for binary cubics and stratify $\PP S_3=\PP^3$.

\begin{proof}[Solution]  Let $X \subseteq \PP^3$ be the rational normal curve.
If $p$ is a point of $\PP^3$ which does not lie on $X$, then considering the image of $X$ in $\PP^2$ by projecting from $p$ we see that either $p$ lies on a tangent line to $X$ or that $p$ lies on a secant line to $X$.  Suppose that $q$ is not a point of $X$ but that $q$ lies on a tangent line.  If $X$ is the twisted cubic and $p= [L^3]$ then $$T_p(X) =\langle [L^3], [L^2M]\rangle $$ where $M$ is any linear form which is not a scalar multiple of $L$. Thus, without loss of generality, to show that $p$ can be written as a sum of $3$ cubes it suffices to show that $x^2y$ is a sum of three cubes.  For this, observe that $$x^2y=\frac{1}{6}((x+y)^3+(y-x)^3-2y^3).$$  Thus, $\PP^3$ is stratified by Waring rank.  Those points of rank $1$ correspond to points of $X$.  Those points of Waring rank $2$ correspond to points which lie on no tangent line to $X$.  Those points of Waring rank $3$ correspond to points which lie on a tangent line and are not on $X$.  All three of these sets are locally closed.
\end{proof}

\subsection*{Exercise \ref{exTWO2}}  Prove the general statement.  If $\Lambda$ is a hyperplane containing $T_{[L^d]}(\nu_d(\PP^n))$ then $\nu_d^{-1}(\Lambda)$ is a degree $d$ hypersurface singular at the point $[L]\in \PP^n$.
\begin{proof}[Solution]
See the solution to Exercise \ref{exONE5}.
\end{proof}

\subsection*{Exercise \ref{exTWO3}}  Solve the big Waring problem for $n=1$ using the double points interpretation.
\begin{proof}[Solution]   We have $N=\binom{d+1}{d}-1=d$.  On the other hand if $p_1,\dots, p_s$ are general points of $\PP^1$ then $\dim (\mathfrak{p}_1^2\cap \dots \cap \mathfrak{p}_s^2)_d = -2s+d+1$ for $1\leq s \leq d+1$.  So $N-2s+d+1\geq N$ implies that $s\geq \lceil \frac{d+1}{2} \rceil$.  In other words $g(1,d)= \lceil \frac{d+1}{2} \rceil$.
\end{proof}

\subsection*{Exercise \ref{exTWO4}}  Show that $\sigma_5(\nu_4(\PP^2))$ is a hypersurface,  \emph{i.e.} that it has dimension equal $13$.
\begin{proof}[Solution] Since $\sigma_5(\nu_4(\PP^2))\not = \PP^{13}$ we conclude $\dim \sigma_5(\nu_4(\PP^2))\geq \dim \sigma_4(\nu_4(\PP^2))+2$.  On the other hand we know $\dim \sigma_5(\nu_4(\PP^2))\leq 13$.  Hence it suffices to show $\dim \sigma_4(\nu_4(\PP^2))\geq 11$.  To see this we note that it costs at most $3$ linear conditions for a plane curve to be singular at a point  hence $\dim(\mathfrak{p}^2_1\cap \dots \mathfrak{p}^2_4)_4\geq 15-12=3$ for all collections of $4$ points.  Applying the double point lemma for a general collection of points, we deduce that $\dim \sigma_4(\nu_4(\PP^2))\geq 11$.
\end{proof}

\subsection*{Exercise \ref{exTWO5}}  Explain the exceptional cases $d=2$ and any $n$.
\begin{proof}[Solution] Let us prove that $g(n,2)=n+1$.  A general $n+1\times n+1$ symmetric matrix has rank $n+1$ and is hence a linear combination of $n+1$  symmetric matrices of rank $1$ and is not a linear combination of any smaller number of rank $1$ symmetric matrices.  On the other hand, symmetric matrices of rank $1$ are exactly those in the image of the quadratic Veronese map.  This explains the exceptional case $d=2$ for all $n$.   \end{proof}

\subsection*{Exercise \ref{exTWO6}}  Explain the exceptional cases $d=4$ and $n=3,4$.
\begin{proof}[Solution] We want to show that $g(3,4)=10$ and not $9$ and that $g(4,4)=15$ and not $14$.  In both cases a
parameter count shows that we can find quadrics through $9$ points and $14$ points in $\PP^3$ and $\PP^4$ respectively.   Squaring these forms produces a quartic singular at these points.  Applying the double point lemma we conclude that $\sigma_9(\nu_4(\PP^3))$ and $\sigma_{14}(\nu_4(\PP^4))$ fail to fill up the space.  On the other hand a parameter count also shows that $\sigma_{10}(\nu_4(\PP^3))$ and $\sigma_{15}(\nu_4(\PP^4))$ fill up the space.
\end{proof}

\subsection*{Exercise \ref{exTWO7}}  Explain the exceptional case $d=3$ and $n=4$.  (Hint: use Castelnuovo's Theorem which asserts that there exists a (unique) rational normal curve passing through $n+3$ generic points in $\PP^n$.)
\begin{proof}[Solution]

By the double point lemma it suffices to pick $7$ general points $p_1,\dots, p_7 \in \PP^4$ and prove that there exists a degree $3$ hypersurface in $\PP^4$ singular at $p_1,\dots, p_7$.  So choose $p_1,\dots, p_7$ general points.  By Castelnuovo's Theorem, there exists a rational quartic curve $X$ in $\PP^4$ passing through the points $p_1,\dots, p_7$.  It suffices to prove that $\sigma_2(X)$ is a degree $3$ hypersurface singular along $X$.

To show that $\sigma_2(X)$ is singular along $X$ if $x\in X$, then we can show that $X\subseteq T_x \sigma_2(X)$ so that $\langle X \rangle \subseteq T_x \sigma_2(X)$.  Since $X$ is non-degenerate we conclude that $T_x \sigma_2(X)=\PP^4$, for all $x\in X$.  Since $\sigma_2(X)$ has dimension $3$ we conclude that $\sigma_2(X)$ is singular along $X$.

To compute the degree we project to $\PP^2$ from a general secant line and count the number of nodes.  Since the resulting curve is rational the number of nodes equals the arithmetic genus of a plane curve of degree $4$ which is $3$.  \end{proof}


\subsection*{Exercise \ref{exTHREE1}} Let $0\not = F\in S_d$. Show that $F^\perp\subset T$ is an ideal and that it is also Artinian. \emph{i.e.}, $(T/F^\perp)_i=0$ for all $i> d$.
\begin{proof}[Solution] If $\partial \in F^{\perp}$ and $\tilde{\partial}\in T$ then by definition of the $T$-action $\tilde{\partial}\partial F= (\tilde{\partial})\partial F=0$  so $\tilde{\partial}\partial \in F^\perp$.   On the other hand since $F$ has degree $d$, $F^\perp$ contains all differential operators of degree greater than $d$.  Hence $T/F^\perp$ is a finite dimensional vector space and hence Artinian.
\end{proof}

\subsection*{Exercise \ref{exTHREE2}}  Show that the map $S_i\times T_i \rightarrow \CC$, $(F,\partial)\mapsto \partial F$ is a perfect paring and that $A=T/F^\perp$ is Artinian and Gorenstein with socle degree $d$.
\begin{proof}[Solution]
 We use the standard monomial basis for $S_i$ and $T_i$.  By definition of the action we have
$$y_0^{a_0}\dots y_n^{a_n}\circ x_0^{b_0}\dots x_n^{b_n} =  \begin{cases}  a_0!\dots a_n! & \text{ iff $a_j=b_j$, for $j=0\dots, n$} \\ 0 & \text{ otherwise }\end{cases} $$ from which the first assertion is clear.

It remains to show that $A$ is Gorenstein with socle degree $d$.  Let $\mathfrak{m}$ denote the homogeneous maximal ideal of $A$.  Then $Soc(A)=\{x \in A \mid x \mathfrak{m} =0\}$.  Using this description we can check that $$\dim Soc(A)_i = \begin{cases}  0 & \textrm{ if $i<d$} \\ 1 & \textrm{ if $i=d$. }\end{cases}$$  Hence $Soc(A)$ is $1$ dimensional and nonzero only in degree $d$ which implies that $A$ is Gorenstein with socle degree $d$.
\end{proof}

\subsection*{Exercise \ref{exTHREE3}}  Given $F\in S_d$ show that $HF(T/F^\perp,t)$ is a symmetric function of $t$.
\begin{proof}[Solution]
Let $A=T/F^\perp$ and let $d$ be the socle degree of $A$. It suffices to show that for $0\leq l\leq d$, multiplication in $A$ defines a perfect pairing $A_{d-l}\times A_l\rightarrow A_d$.  If $l=d$ or $l=0$ the assertion is clear.  Let $0<l<d$.  Let $y\in A_l$.  Suppose $yx=0$ for all $x\in A_{d-l}$.  Let's prove that $y\in Soc(A)$.  Since $Soc(A)$ is zero in degrees less than $d$ we will arrive at a contradiction.

To prove that $y\in Soc(A)$ it suffices to prove that $y$ annihilates every homogeneous element of $A$ which has positive degree.  If $M\in A_n$ and $n>d-l$ then $yM \in A_{l+n}=0$ since $l+n>d$.  On the other hand we have by assumption that $yA_{d-l}=0$.  Descending induction with base case $d-l$ proves that $yA_l=0$ for $0<l<d-l$.  Indeed suppose $0<n<d-l$ and let $M\in A_n$.  Suppose that $yM\not = 0$.  We have $Mx_i\in A_{n+1}$ for $i=0,\dots, n$.  By induction $yMx_i=0$.  Hence $yM\mathfrak{m}=0$ so $yM\in Soc(A)$.  Since $\deg yM <d$ this is a contradiction.  Hence $y\in Soc(A)$ which is also a contradiction.
\end{proof}

\subsection*{Exercise \ref{exTHREE4}} Use the Apolarity Lemma to compute $\operatorname{rk}(x_0x_1^2)$.  Then try the binary forms $x_0x_1^d$.
\begin{proof}[Solution]  $F^\perp = \langle \partial_{x_0}^2, \partial_{x_1}^{d+1} \rangle$ so $\operatorname{rank} F = d+1$.
\end{proof}
\subsection*{Exercise \ref{exTHREE5}}  Use the Apolarity Lemma to explain the Alexander-Hirschowitz exceptional cases.
\begin{proof}[Solution]  We explain the exceptional cases $d = 4$ and $n = 2,3$ or $4$.  Then exceptional case $d = 3$ and $n = 4$ can be treated via syzygies.

Since $\binom{2 + 2}{2} = 6 > 5$ if $I$ is the ideal of $5$ points in $\PP^2$, then $I$ contains a quadric.
Since $\binom{2 + 3}{ 2} = 10 > 9$ if $I$ is the ideal of $9$ points in $\PP^3$, then $I$ contains a quadric.  Since $\binom{2 + 4}{2} = 15 > 14$ if $I$ is the ideal of $14$ points in $\PP^4$, then $I$ contains a quadric.

On the other hand, if $F$ is a general form of degree $4$, in $\CC[x_0,..,x_n]$, for $n = 2,3$ or $4$, then $F^\perp$ contains no quadrics.  We work out the case $n = 2$ explicitly.   The case $n = 3$ or $4$ is similar.

Let $S = \CC[x,y,z]$.  Every form $F \in S_4$ determines a linear map from the vector space of differential operators of degree $2$ to the space of degree $2$ polynomials in $S$.  This map is determined by applying the differential operators to $F$.
Explicitly, if
\begin{multline*}
F = a x^4 + b x^3 y + c x^3 z + d x^2 y^2 + e x^2 y z + f x^2 z^2 + \\
  g x y^3 + h x y^2 z + i x y z^2 + j x z^3 + k y^4 + l y^3 z + \\
  m y^2 z^2 + o y z^3 + p z^4
\end{multline*}
then, using the basis
$$\partial_{xx}, \partial_{xy}, \partial_{xz}, \partial_{y^2}, \partial_{yz}, \partial_{z^2}$$ for the source, and the basis
$$x^2, xy, xz, y^2, yz, z^2$$ for the target space, the matrix for this map is given by

$$
\left[
\begin{array}{cccccc}
12 a & 3 b & 3 c & 2 d & e & 2 f \\
6 b & 4 d & 2 e & 6 g & 2 h & 2 i \\
6 c & 2 e & 4 f & 2 h & 2 i & 6 j \\
2 d & 3 g & h & 12 k & 3 l & 2 m \\
2 e & 2 h & 2 i & 6 l & 4 m & 6 o \\
2 f & i & 3 j & 2 m & 3 o & 12 p
\end{array} \right] \text{.}
$$
Elements in the kernel of this map correspond to elements of $(F^\perp)_2$.
The collection of forms for which this map is injective is given by the non-vanishing of the determinant of this matrix.  We conclude that $(F^\perp)_2$ is zero for a general quartic form in $S_4$.

Note that we have shown that if $F$ is a general form of degree $4$ in $ \CC[x_0,..,x_n]$, for $n = 2,3$ or $4$, then

\[
\begin{array}{llllllll}

t               & 0 &   1   &   2   &   3   & 4  & 5 \\

\hline
\\
HF(T/F^\perp,t) & 1 & n + 1 &  \binom{n + 2}{2}   &  n + 1&  1 & 0 & \rightarrow
\end{array}
\]

On the other hand, for every ideal $I$ of respectively, $5$ points in $\PP^2$, $9$ points in $\PP^3$,  or $14$ points in $\PP^4$, we have shown that
$$HF(T/ I, 2) < \binom{n + 2}{2} \text{.}$$

Thus, using the Apolarity Lemma, a general quartic form in $\CC[x_0,..,x_n]$, for $n = 2,3$ or $4$, cannot have rank respectively, $5$, $9$, or $14$.

To see an issue which arises in the case exceptional case $d = 3$ and $n = 4$, note that if $F$ is any cubic form in $\CC[x_0,x_1,x_2,x_3,x_4]$ then
\[
\begin{array}{llllll}

t               & 0 &   1   &   2   &   3   & 4  \\

\hline
\\
HF(T/F^\perp,t) & 1 & 5 &  5 & 1 &  0  \rightarrow
\end{array}
\]
Thus the Hilbert function  $HF(T/F^\perp,t)$ is the same for all $$F \in \CC[x_0,x_1,x_2,x_3,x_4]_3 \text{.}$$
Nevertheless, the exceptional case $d = 3$ and $n = 4$ can still be explained using the Apolarity Lemma, although a more detailed study is needed to conclude that $F^\perp$ for a general $F \in \CC[x_0,x_1,x_2,x_3,x_4]_3 $ contains no ideal of $7$ points in $\PP^3$.  See the paper \cite{RanestadSchreyer_VSP} for a more detailed discussion.
\end{proof}

\subsection*{Exercise \ref{exTHREE6}} Compute $\operatorname{rk}(F)$ when $F$ is a quadratic form.
\begin{proof}[Solution]  $\operatorname{rk}(F)= \operatorname{rank} M_F$ where $M_F$ is the symmetric matrix associated to $F$. (See Exercise \ref{exONE2+++}.)
\end{proof}

\subsection*{Exercise \ref{exTHREE7}} Prove that $\operatorname{rk}(L^d+M^d+N^d)=3$ whenever $L,M$ and $N$ are linearly independent linear forms.
\begin{proof}[Solution] It is clear that $\operatorname{rk}(L^d+M^d+N^d)\leq3$.  To show that $\operatorname{rk}(L^d+M^d+N^d)$ is not less than $3$,  without loss of generality it suffices to consider the case that $F=x_0^d+x_1^d+x_2^d$.  In this case
  $F^{\perp}=\langle \partial_{x_0}^{d+1},\partial_{x_1}^{d+1}, \partial_{x_2}^{d+1}, \partial_{x_3}, \dots, \partial_{x_n}  \rangle $.  Now if $F^{\perp}$ contained the ideal of $1$ or $2$ distinct points in $\PP^n$ then $\langle \partial_{x_0}^{d+1},\partial_{x_1}^{d+1}, \partial_{x_2}^{d+1}\rangle$ would contain a linear form which it does not.
\end{proof}


\subsection*{Exercise \ref{exFOUR1}} Workout a matrix representation for the Segre
varieties with two factors
$\mathbb{P}^{n_1}\times\mathbb{P}^{n_2}$.
\begin{proof}[Solution]
Let $A\cong \CC^{n_{1}+1}$ and $B\cong \CC^{n_{2}+1}$. By
definition, points in
$\Seg(\mathbb{P}^{n_1}\times\mathbb{P}^{n_2})$ are of the form
\[[a\otimes b] \in \PP(A\otimes B)
.\]

Now consider $A\otimes B$ as a space of matrices $A^{*}\to B$. Then, by choosing bases of $A$ and $B$, we may represent $a$ as a column vector $(a_{1},\dots,a_{n_{1}+1})^{t}$ (an element of $A^{*}$) and $b$ as a row vector $(b_{1},\dots,b_{n_{2}+1})$ so that
the tensor product $a\otimes b$ becomes the product of a column and a row:
\[
a\otimes b = (a_{1},\dots,a_{n_{1}+1})^{t} \cdot (b_{1},\dots,b_{n_{2}+1})
=
\begin{pmatrix}a_{1}b_{1} & a_{1}b_{2} & \dots & a_{1}b_{n_{2}+1} \\
a_{2}b_{1} & a_{2}b_{2} & \dots &  a_{2}b_{n_{2}+1} \\
\vdots & \ddots & & \vdots \\
a_{n_{1} + 1}b_{1} & a_{n_{1} + 1}b_{2} & \dots &  a_{n_{1} + 1}b_{n_{2}+1}
\end{pmatrix}.
\]
So we see that (up to scale) elements of $\Seg(\PP^{n_{1}}\times \PP^{n_{2}})$ correspond to rank-one $(n_{1}+1) \times (n_{2}+1)$ matrices.
\end{proof}

\subsection*{Exercise \ref{exFOUR2}} Let $X=\PP V_{1}\times \dots \times \PP V_{t}$ and let $[v] =[v_{1}\otimes\dots\otimes v_{t}]$ be a point of $X$.
Show that the cone over the tangent space to $X$ at $v$ is the
span of the following vector spaces:
\[\begin{matrix}
V_{1}\otimes v_{2}\otimes v_{3}\otimes \dots \otimes v_{t},\\
v_{1}\otimes V_{2}\otimes v_{3}\otimes \dots \otimes v_{t},\\
\vdots \\
v_{1}\otimes v_{2}\otimes \dots \otimes v_{t-1}\otimes V_{t}.
\end{matrix}\]
\begin{proof}[Solution]
The cone over the tangent space to a variety $X$ at a point $v$
may be computed by considering all curves $\gamma:[0,1]\to
\widehat X$ such that $\gamma(0)=v$, and taking the linear span of
all derivatives at the origin:
\[
\widehat{T_{x}X} = \left\{
\gamma'(0) \mid \gamma:[0,1]\to \widehat X, \gamma(0)=v
\right\}
.\]

Now take $X=\PP V_{1}\times \dots \times \PP V_{t}$ and let $[v] =[v_{1}\otimes\dots\otimes v_{t}]$ be a point of $X$.

All curves $\gamma(t)$ on $\widehat X$ through $v$ are of the form
$\gamma(t) =v_{1}(t)\otimes\dots\otimes v_{t}(t)$, where $v_{i}(t)$ are curves in $V_{i}$ such that $v_{i}(0)=v_{i}$.

Now apply the product rule, and for notational convenience, set $v'_{i} = v_{i}'(0) \in V_{i}$.  We have
\[
\gamma'(0) = v_{1}'\otimes v_{2}\otimes\dots\otimes v_{t}
+
v_{1}\otimes v_{2}'\otimes v_{3}\otimes\dots\otimes v_{t}
+\dots +
v_{1}\otimes\dots\otimes v_{t}'
.\]

Since $v_{i}'$ can be anything in $V_{i}$ we get the result.
\end{proof}

\subsection*{Exercise \ref{exFOUR3}} Show that $\sigma_{2}(\PP^{1}\times \PP^{1}\times \PP^{1}) = \PP^{7}$.
\begin{proof}[Solution]
Let $a\otimes b\otimes c + \tilde a\otimes\tilde  b \otimes\tilde
c$ be a general point on $\sigma_{2}(\PP A\times \PP B \times \PP
C)$, with $A\cong B \cong C \cong \CC^{2}$.

By the previous exercise, we have
\[
\widehat{T_{a\otimes b\otimes c}}(\PP^{1}\times \PP^{1}\times \PP^{1}) =
A\otimes b\otimes c  + a\otimes B \otimes c + a \otimes b \otimes C
,\]
and similarly
\[
\widehat{T_{\tilde a\otimes\tilde  b\otimes\tilde  c}}(\PP^{1}\times \PP^{1}\times \PP^{1}) =
A\otimes \tilde b\otimes\tilde  c  +\tilde  a\otimes B \otimes\tilde c +\tilde a \otimes\tilde b \otimes C
.\]

Now by Terracini's lemma we have
\begin{multline}\label{eq:tan}
\widehat{T_{a\otimes b\otimes c + \tilde a\otimes\tilde  b \otimes\tilde  c}}\sigma_{2}(\PP^{1}\times \PP^{1}\times \PP^{1})\\
=
A\otimes b\otimes c  + a\otimes B \otimes c + a \otimes b \otimes C
+
A\otimes \tilde b\otimes\tilde  c  +\tilde  a\otimes B \otimes\tilde c +\tilde a \otimes\tilde b \otimes C
.\end{multline}
Because we chose a general point, we have $\{a,\tilde b\} = A$ and similarly for $B$ and $C$.
Now consider the linear space in \eqref{eq:tan}. We see that we can get every tensor monomial in $A\otimes B\otimes C$ -- all monomials with 0 or 1 $\tilde{}$ occur in the first 3 summands, while all monomials with 2 or 3 $\tilde{}$ occur in the second 3 summands.

So the cone over the tangent space at a general point is 8 dimensional, so the secant variety (being irreducible) fills the whole ambient space.
\end{proof}

\subsection*{Exercise \ref{exFOUR4}}Use the above description of the tangent space of the Segre product and Terracini's lemma to show that $\sigma_{3}(\PP^{1}\times \PP^{1}\times \PP^{1}\times \PP^{1})$ is a hypersurface in $\PP^{15}$ and not the entire ambient space as expected. This shows that the four-factor Segre product of $\PP^{1}$s is defective.
\begin{proof}[Solution]
 This should be done on the computer.  Take the following to be a general point:
\begin{multline*}
a\otimes b\otimes c\otimes d +
\tilde a\otimes\tilde b\otimes\tilde c\otimes\tilde d +\\
(\alpha_{1}a+\alpha_{2}\tilde a)\otimes(\beta_{1}b+\beta_{2}\tilde b)\otimes(\gamma_{1}c+\gamma_{2}\tilde c)\otimes(\delta_{1}d+\delta_{2}\tilde d)
,\end{multline*}
where $[\alpha_{1},\alpha_{2}],[\beta_{1},\beta_{2}],[\gamma_{1},\gamma_{2}],[\delta_{1},\delta_{2}] \in \PP^{1}$.

Then compute the derivatives at the origin.  Keep track of all the parameters and letting the $a(0)'$ $\tilde{a}(0)'$ etc., vary produce 16 vectors that span the tangent space.  Now compute the rank of the matrix with these vectors as columns.  It gets hard to do by hand, but the next section produces an easier way to do the problem.

The ``easier way'' is to consider all 3 essentially different 2-flattenings, and show that two of them are algebraically independent.
\end{proof}

\subsection*{Exercise \ref{exFOUR5}}
Show that $T$ has rank 1 if and only if its multilinear rank is
$(1,\dots,1)$.
\begin{proof}[Solution]
 It is equivalent (by taking transposes) to consider the $n-1$ flattenings.
Let $\varphi_{i,T}:(V_{1}\otimes \dots\otimes \widehat{V_{i}} \otimes \dots\otimes V_{n})^{*} \to V_{i}$
denote the $(n-1)$-flattening to the $i$th factor.  Denote by $A_{i}$ the image of $\varphi_{i,T}$.

If $\varphi_{i,T}$ has rank 1 then $\dim A_{i} = 1$. So, we must have $T \in A_{1}\otimes \dots A_{n}$.  But every tensor in $A_{1}\otimes A_{2}\otimes \dots A_{n}$ has rank 1 if all the factors have dimension 1.

Conversely, if $T=a_{1}\otimes \dots a_{n} \in V_{1}\otimes \dots \otimes V_{n}$ the image of $\phi_{i,T}$  is the line through $a_{i}$, so the multilinear rank is $(1,1,\dots,1)$.
\end{proof}

\subsection*{Exercise \ref{exFOUR6}} Let $X=\PP V_{1}\times \dots \times \PP V_{t}$. Show that if $r
\leq r_{i}$ for $1\leq i \leq t$,
\[\sigma_{r}(X)\subset \Sub_{r_{1},\dots,r_{t}}.\]
\begin{proof}[Solution]
 A general point on $\sigma_{r}(X)$ is $p=\sum_{s=1}^{r} \bigotimes_{i=1}^{t} a_{i,s}$,
where for fixed $i$ the $a_{i,s} \in V_{i}$ are linearly independent. Set $A_{i} = \{
a_{i,1},\dots,a_{i,r}\}$.  Then $p \in A_{1}\otimes \dots \otimes A_{t}$, so $p\in \Sub_{r_{1},\dots,r_{t}}$.  Now take the orbit closure of $p$ to obtain the result.
\end{proof}

\subsection*{Exercise \ref{exFOUR7}} \begin{enumerate}
\item Show that if $T$ has rank $1$ then $\varphi_{T}$ has rank 2.
\item Show that $\varphi$ is additive in its argument, \emph{i.e.}
show that $\varphi_{T+T'} = \varphi_{T}+ \varphi_{T'}.$
\end{enumerate}
\begin{proof}[Solution]
 If $T$ has rank 1, after change of coordinates we may assume that $T = v_{1}\otimes w_{1}\otimes x_{1}$. Then the matrix $\varphi_{T}$ has precisely two ones in different rows and columns, so clearly has rank 2.

For the second part, notice
\[\varphi_{T} + \varphi_{T'} =
\begin{pmatrix}
0 & T^{1} & -T^{2} \\
-T^{1} & 0 &  T^{3} \\
T^{2} & -T^{3} &0
\end{pmatrix}
+
\begin{pmatrix}
0 & T'^{1} & -T'^{2} \\
-T'^{1} & 0 &  T'^{3} \\
T'^{2} & -T'^{3} &0
\end{pmatrix}
\]
\[=
\begin{pmatrix}
0 & T^{1} + T'^{1} & -T^{2} + T'^{2} \\
-T^{1} + T'^{1} & 0 &  T^{3} + T'^{3} \\
T^{2} + T'^{2} & -T^{3} + T'^{3} &0
\end{pmatrix}
=
\varphi_{T+T'}.
\]
\end{proof}


\section{Looking forward, further readings}

This series of lectures draws from many sources from a large group of authors. We do our best to collect representative works here so that the reader may have some starting points for further study. A general introductory reference for the material we treated is the booklet is \cite{MR1381732}.  Here we collect a few guiding questions and an extensive list of references.

\subsection*{Guiding questions}
In our opinion, there are a few leading topics that are still driving current research in the area of secant varieties. These topics are dimension, identifiability, decomposition, and equations. More specifically here are 4 leading questions:
\begin{enumerate}
\item For a variety  $X$ what are the dimensions of the higher secant varieties to $X$? There is much interest when $X$ is the variety of elementary tensors of a given format (partially symmetric, skew symmetric, general).
\item Suppose $p\in \sigma_{s}(X)$ is general, when does $p$ have a \emph{unique} representation as a sum of $s$ points from $X$? When uniqueness occurs we say that $\sigma_{s}(X)$ is \emph{generically identifiable}.   
\item Suppose $X \subset \CC^{N}$ and $p\in \CC^{N}$. If $X$ is not degenerate, then we know that there is some $s$ so that $p$ has a representation as the linear combination of $s$ points from $X$. For special $X$ (Segre, Veronese, etc.,) find: (a) determine the minimal $s$ explicitly, and (b)find efficient algorithms when $s$ is relatively small to find such a decomposition of $p$. 
\item How do we find equations $\sigma_{s}(X)$ in general? What is the degree of $\sigma_{s}(X)$? Is $\sigma_{s}(X)$ a Cohen-Macaulay variety?  Again, there is much interest when $X$ is the variety of elementary tensors of a given format.
\end{enumerate}

\subsection*{Background material}

Textbooks:
\cite{GKZ,bruns-vetter,Burgisser, CGHRR,CLO_text, Eisenbud_syzygies, Eisenbud_tome, Enriques_book,FlennerOCarrol_Joins, FultonHarris, gantmacher, GoodWall, Harris, Hartshorne_tome, LandsbergTensorBook,OSS_text,Pachter-Sturmfels,Weyl,Weyman,Zak}.

Classical: \cite{Terracini}. See also the nice overview in
\cite{MR1381732} and the introduction of
\cite{RanestadSchreyer_VSP}.

\subsection*{ Dimensions of secant varieties}

For Veronese varieties, the capstone result is that of Alexander and Hirschowitz, \cite{AH95}. Ottaviani and Brambilla \cite{OttBra08_AH} provided a very nice exposition.
Related work on polynomial interpolation: \cite{AH92,OttBra11_interpolation}. 

Waring's problem for binary forms: 
\cite{IarrobinoKanev_text, ComSei11_binary, LandsbergTeitler_SymRank} 
Waring's problem for polynomials: 
\cite{
FrobergOttavianiShapiro_Waring, Mella06_Waring, Mella09_Waring, Ottaviani09_Waring, RanestadSchreyer_VSP, Comon-Mourrain_Waring, Chipalkatti04a_Waring}.
Waring's problem for monomials was solved in \cite{CCG_Waring} and
an alternative proof can be found also in \cite{BBT_Waring} where
the apolar sets of points to monomials are described. See also
\cite{BCG_monomials}.

The polynomial Waring problem over the reals was investigated by
Comon and Ottaviani, \cite{ComonOtt_Binary}, with some solutions
to their questions provided by \cite{Blekherman12_RealRank,
Causa_rank_binary, BallicoTypical}. The case of real monomials in
two variables is discussed in \cite{BCG_monomials}. For typical
ranks of tensors see \cite{ ComonTenBerg_Typical,
Friedland_Typical,SumiSakata_Typical}.

Recent algorithms for Waring decomposition: \cite{IarrobinoKanev_text, OedOtt13_Waring,BGI11_Waring,BCMT_Waring,CGLM_Waring}.

The notion of Weak Defectivity:
\cite{ChiCil_WeakDefective,ChiCil06_ksecant,Ballico2005_weak}.

Dimensions of Secant varieties of Segre-Veronese varieties: \cite{CGG2_Segre, CGG2.5, AOP_Segre, CGG_P1s, CGG4_Segre, CGG5_rational, abresica_SegreVeronese, CGG3_SegreVeronese,CGG1_tensors, AboBram09, AboBram12, AboBram13, BCC_PmPn, BBC_PnP1}.
Nice Summary of results on dimensions of secant varieties:  \cite{Catalisano_slides}.

\subsection*{ Equations of secant varieties} See \cite{Ott07_Luroth,LanMan04_Seg,
LanWey_Seg, RaicuGSS, SidmanSullivant_Prolongations, LanWey_secant_CHSS, LanOtt11_Equations, RaicuCat, Strassen83_rank, CEO_partial,LanMan08_Strassen,
Bernardi_Ideal_Sym,LanWey_Tan, OedingTan, Kanev_catalecticant,
RaicuGSS, Raicu_thesis, SidmanVermeire_Secant}.

The Salmon Problem:
\cite{prize,Friedland2010_salmon,OedingBates,Friedland-Gross2011_salmon}.

\subsection*{ Applications}
The question of best low-rank approximation of tensors is often ill posed \cite{deSilva_Lim_ill_posed}, and most tensor problems are NP-Hard, \cite{HillarLim09mosttensor}.

Algebraic Statistics: \cite{Drton-Sturmfels-Sullivant, DMS_book, GSS, SturmfelsOpenProblems}.

Phylogenetics: \cite{AllmanRhodes03,AllmanRhodes08, Casanellas_phylo, ERSS_phylo}.

Signal Processing: see
\cite{AFCC_blind,Comon1_canonical,Comon2_tensor,ComoR06_blind,DD08_blind,DFDV02_block_decomp,LathCastCard_cumulant,LathCast_blind,lathauwer_Schur,
lathauwer:642_matrix_diag,DBLP:conf/coco/2007}.


Matrix Multiplication:
\cite{Lan_2by2_matrix,Lan_matrix_mult,Landsberg_explicit,Landsberg_NewLower,LanOtt_NewLower}.

\subsection*{ Related varieties}
Grassmannians: \cite{ CGG6_Grassman,AOP_Grassman,BorBuc_Lagrangian}.

Discriminants and Hyperdeterminants are intimately related to secant varieties of Segre-Veronese varieties, see \cite{CCDRS_Discriminants,GKZ89,GKZ92,HSYY_hyperdet,Oeding_Hyperdet,WeymanBoffi_Hyperdet,WeymanZelevinsky_mult, WeymanZelevinsky_sing}.

Chow varieties and monomials\cite{ArrondoBernardi_chow,Briand_chow,Carlini_chow}.

\subsection*{ Related concepts}

Eigenvectors of Tensors:
\cite{Lim05_evectors,Qi05_eigen,BKP_tensor_evals,Hu_Edet,OttSturm_Eigenvectors,NiQiWang_Epoly,CartwrightSturmfels2011,KoldaMayo_tensor_eigenpairs}.

Veronese reembeddings: \cite{BucBuc_Veronese}.
Hilbert Schemes: \cite{ErmanVelasco_smoothable}.

Orbits: \cite{Djokovic,VinbergElashvili_classification}.

Ranks and Decompositions:
\cite{BallicoBernardi2010_decomp,BucLan_ranks,
BucLan_3rdSecant}.

Asymptotic questions: \cite{draisma-kuttler_bounded}.

\subsection*{ Software}
Symbolic computation: \cite{Singular,CoCoA, M2, lie}.

Numerical Algebraic Geometry: \cite{HOM4PS,PHC,regen,BHSW_Bertini,Bertini,SomeseWampler05,Li_homotopy}.

\bibliographystyle{alpha}

\begin{thebibliography}{DLFDMV02}

\bibitem[AB09]{AboBram09}
H. Abo and M.~C. Brambilla.
\newblock Secant varieties of {S}egre-{V}eronese varieties {$\mathbb{P}^m\times
  \mathbb{P}^n$} embedded by {$\mathcal{O}(1,2)$}.
\newblock {\em Experiment. Math.}, 18(3):369--384, 2009.


\bibitem[AB12]{AboBram12}
\bysame. 
\newblock New examples of defective secant varieties of {S}egre-{V}eronese
  varieties.
\newblock {\em Collect. Math.}, 63(3):287--297, 2012.

\bibitem[AB13]{AboBram13}
\bysame. 
\newblock On the dimensions of secant varieties of {S}egre-{V}eronese
  varieties.
\newblock {\em Ann. Mat. Pura Appl. (4)}, 192(1):61--92, 2013.

\bibitem[AOP09]{AOP_Segre}
H.~Abo, G.~Ottaviani, and C.~Peterson.
\newblock Induction for secant varieties of {S}egre varieties.
\newblock {\em Trans. Amer. Math. Soc.}, 361(2):767--792, 2009.

\bibitem[AOP12]{AOP_Grassman}
\bysame. 
\newblock Non-defectivity of {G}rassmannians of planes.
\newblock {\em J. Algebraic Geom.}, 21(1):1--20, 2012.

\bibitem[Abr08]{abresica_SegreVeronese}
S. Abrescia.
\newblock About the defectivity of certain {S}egre-{V}eronese varieties.
\newblock {\em Canad. J. Math.}, 60(5):961--974, 2008.

\bibitem[AFCC04]{AFCC_blind}
L.~Albera, A.~Ferreol, P.~Comon, and P.~Chevalier.
\newblock Blind identification of overcomplete mixtures of sources ({BIOME}).
\newblock {\em Lin. Algebra Appl.}, 391:3--30, November 2004.

\bibitem[AH92]{AH92}
J.~Alexander and A.~Hirschowitz.
\newblock La m\'ethode d'{H}orace \'eclat\'ee: application \`a l'interpolation
  en degr\'e quatre.
\newblock {\em Invent. Math.}, 107(3):585--602, 1992.

\bibitem[AH95]{AH95}
\bysame. 
\newblock Polynomial interpolation in several variables.
\newblock {\em J. Algebraic Geom.}, 4(2):201--222, 1995.

\bibitem[All]{prize}
E.~Allman.
\newblock  \emph{Open Problem: Determine the ideal defining
  {$Sec_{4}(\mathbb{P}^{3} \times \mathbb{P}^{3} \times \mathbb{P}^{3})$}},
  \url{http://www.dms.uaf.edu/~eallman/salmonPrize.pdf},
  2010.


\bibitem[AR03]{AllmanRhodes03}
E.~Allman and J.~Rhodes.
\newblock Phylogenetic invariants for the general {M}arkov model of sequence
  mutation.
\newblock {\em Math. Biosci.}, 186(2):113--144, 2003.

\bibitem[AR08]{AllmanRhodes08}
\bysame. 
\newblock Phylogenetic ideals and varieties for the general {M}arkov model.
\newblock {\em Adv. in Appl. Math.}, 40(2):127--148, 2008.

\bibitem[AB11]{ArrondoBernardi_chow}
E. Arrondo and A. Bernardi.
\newblock On the variety parameterizing completely decomposable polynomials.
\newblock {\em Journal of Pure and Applied Algebra}, 215(3):201--220, 2011.

\bibitem[BKP11]{BKP_tensor_evals}
G.~Ballard, T.~Kolda, and T.~Plantenga.
\newblock Efficiently computing tensor eigenvalues on a {GPU}, 2011.
\newblock CSRI Summer Proceedings.

\bibitem[Bal05]{Ballico2005_weak}
E.~Ballico.
\newblock On the weak non-defectivity of {V}eronese embeddings of projective
  spaces.
\newblock {\em Cent. Eur. J. Math.}, 3(2):183--187 (electronic), 2005.

\bibitem[{Bal}12]{BallicoTypical}
\bysame. 
\newblock {On the typical rank of real bivariate polynomials}.
\newblock {\em ArXiv e-prints}, April 2012. 


\bibitem[BB11a]{BallicoBernardi2010_decomp}
E.~Ballico and A.~Bernardi.
\newblock Decomposition of homogeneous polynomials with low rank.
\newblock {\em Mathematische Zeitschrift}, 271(3-4):1141--1149, 2011.

\bibitem[BBC12]{BBC_PnP1}
E. Ballico, A. Bernardi, and M. V. Catalisano.
\newblock Higher secant varieties of {$\Bbb P^n\times\Bbb P^1$} embedded in
  bi-degree {$(a,b)$}.
\newblock {\em Comm. Algebra}, 40(10):3822--3840, 2012.


\bibitem[BHSW08]{BHSW_Bertini}
D.J. Bates, J.D. Hauenstein, A.J. Sommese, and C.W. Wampler.
\newblock Software for numerical algebraic geometry: a paradigm and progress
  towards its implementation.
\newblock In {\em Software for algebraic geometry}, volume 148 of {\em IMA Vol.
  Math. Appl.}, pages 1--14. Springer, New York, 2008.


\bibitem[BHSW]{Bertini}
\bysame. 
Bertini: Software for Numerical Algebraic Geometry.
Available at \url{http://www.nd.edu/~sommese/bertini}, 2010.
  
\bibitem[BO11a]{OedingBates}
D.~J. Bates and L. Oeding.
\newblock Toward a salmon conjecture.
\newblock {\em Exp. Math.}, 20(3):358--370, 2011.

\bibitem[Ber08]{Bernardi_Ideal_Sym}
A.~Bernardi.
\newblock Ideals of varieties parameterized by certain symmetric tensors.
\newblock {\em J. Pure Appl. Algebra}, 212(6):1542--1559, 2008.

\bibitem[BCC11]{BCC_PmPn}
A. Bernardi, E. Carlini, and M. V. Catalisano.
\newblock Higher secant varieties of {$\Bbb P^n\times\Bbb P^m$} embedded in
  bi-degree {$(1,d)$}.
\newblock {\em J. Pure Appl. Algebra}, 215(12):2853--2858, 2011.


\bibitem[BGI11]{BGI11_Waring}
A.~{Bernardi}, A.~{Gimigliano}, and M.~{Id{\`a}}.
\newblock Computing symmetric rank for symmetric tensors.
\newblock {\em Journal of Symbolic Computation}, 46(1):34--53, 2011.

\bibitem[{Ble}12]{Blekherman12_RealRank}
G.~{Blekherman}.
\newblock {Typical Real Ranks of Binary Forms}.
\newblock {\em ArXiv e-prints}, May 2012.

\bibitem[BW00]{WeymanBoffi_Hyperdet}
G.~Boffi and J.~Weyman.
\newblock Koszul complexes and hyperdeterminants.
\newblock {\em J. Algebra}, 230(1):68--88, 2000.

\bibitem[BCG11]{BCG_monomials}
M. Boij, E. Carlini, and A.~V. Geramita.
\newblock Monomials as sums of powers: the real binary case.
\newblock {\em Proc. Amer. Math. Soc.}, 139(9):3039--3043, 2011.

\bibitem[BB11b]{BorBuc_Lagrangian}
A. Boralevi and J. Buczy{\'n}ski.
\newblock Secants of {L}agrangian {G}rassmannians.
\newblock {\em Ann. Mat. Pura Appl. (4)}, 190(4):725--739, 2011.

\bibitem[BCMT10]{BCMT_Waring}
J. Brachat, P. Comon, B. Mourrain, and E. Tsigaridas.
\newblock Symmetric tensor decomposition.
\newblock {\em Linear Algebra and its Applications}, 433(11-12):1851--1872,
  2010.


\bibitem[BO08]{OttBra08_AH}
M.C.~Brambilla and G.~Ottaviani.
\newblock On the {A}lexander-{H}irschowitz theorem.
\newblock {\em J. Pure Appl. Algebra}, 212(5):1229--1251, 2008.


\bibitem[BO11b]{OttBra11_interpolation}
\bysame. 
\newblock On partial polynomial interpolation.
\newblock {\em Linear Algebra and its Applications}, 435(6):1415 -- 1445, 2011.

\bibitem[Bri10]{Briand_chow}
E.~Briand.
\newblock Covariants vanishing on totally decomposable forms.
\newblock In {\em Liaison, {S}chottky problem and invariant theory}, volume 280
  of {\em Progr. Math.}, pages 237--256. Birkh\"auser Verlag, Basel, 2010.

\bibitem[BV88]{bruns-vetter}
W. Bruns and U. Vetter.
\newblock {\em Determinantal rings}, volume 1327 of {\em Lecture Notes in
  Mathematics}.
\newblock Springer-Verlag, Berlin, 1988.
\bibitem[BB10]{BucBuc_Veronese}
W.~{Buczy{\'n}ska} and J.~{Buczy{\'n}ski}.
\newblock {Secant varieties to high degree Veronese reembeddings, catalecticant
  matrices and smoothable Gorenstein schemes}.
\newblock {\em ArXiv:1012.3563}, December 2010.
\newblock to appear: Journal of Algebraic Geometry.


\bibitem[BBT13]{BBT_Waring}
W. Buczy{\'n}ska, J. Buczy{\'n}ski, and Zach Teitler.
\newblock Waring decompositions of monomials.
\newblock {\em J. Algebra}, 378:45--57, 2013.

\bibitem[BL11]{BucLan_3rdSecant}
J.~{Buczy{\'n}ski} and J.~M. {Landsberg}.
\newblock {On the third secant variety}.
\newblock {\em ArXiv e-prints}, November 2011.
\newblock to appear: Journal of Algebraic Combinatorics.






\bibitem[BL13]{BucLan_ranks}
J. Buczy{\'n}ski and J.M. Landsberg.
\newblock Ranks of tensors and a generalization of secant varieties.
\newblock {\em Linear Algebra Appl.}, 438(2):668--689, 2013.

\bibitem[BCS97]{Burgisser}
P.~B{\"u}rgisser, M.~Clausen, and M.~Shokrollahi.
\newblock {\em Algebraic complexity theory}, volume 315 of {\em Grundlehren der
  Mathematischen Wissenschaften [Fundamental Principles of Mathematical
  Sciences]}.
\newblock Berlin: Springer-Verlag, 1997.


\bibitem[Car05]{Carlini_chow}
E.~Carlini.
\newblock Codimension one decompositions and {C}how varieties.
\newblock In {\em Projective varieties with unexpected properties}, pages
  67--79. Walter de Gruyter GmbH \& Co. KG, Berlin, 2005.


\bibitem[CCG12]{CCG_Waring}
E. Carlini, M.~V. Catalisano, and A.~V. Geramita.
\newblock The solution to the {W}aring problem for monomials and the sum of
  coprime monomials.
\newblock {\em J. Algebra}, 370:5--14, 2012.

\bibitem[CK11]{MR2776439}
E. Carlini and J. Kleppe.
\newblock Ranks derived from multilinear maps.
\newblock {\em J. Pure Appl. Algebra}, 215(8):1999--2004, 2011.

\bibitem[CEO12]{CEO_partial}
D.~{Cartwright}, D.~{Erman}, and L.~{Oeding}.
\newblock {Secant varieties of ${\bf P}^{2} \times {\bf P}^{n}$ embedded by
  ${\mathcal O}(1,2)$}.
\newblock {\em Journal of the London Mathematical Society}, 85(1):121--141,
  2012.

\bibitem[CS11a]{CartwrightSturmfels2011}
D.~Cartwright and B.~Sturmfels.
\newblock The number of eigenvalues of a tensor.
\newblock {\em Linear Algebra and its Applications}, 438(2):942--952, 2013.

\bibitem[CFS11]{Casanellas_phylo}
M. Casanellas and J. Fern{\'a}ndez-S{\'a}nchez.
\newblock Relevant phylogenetic invariants of evolutionary models.
\newblock {\em J. Math. Pures Appl. (9)}, 96(3):207--229, 2011.

\bibitem[Cat]{Catalisano_slides}
M. V. Catalisano. Higher secant varieties of Segre, Segre-Veronese and Grassmann varieties.
\emph{Interactions between Commutative Algebra and Algebraic Geometry
the 22nd annual Route 81 conference --- in honour of Tony Geramita} [Conference], Kingston, Ontario. 20 Oct. 2012.
\url{http://www.mast.queensu.ca/~ggsmith/Route81/catalisano.pdf}.

\bibitem[CGG02a]{CGG1_tensors}
M.~V. Catalisano, A.~V. Geramita, and A.~Gimigliano.
\newblock On the rank of tensors, via secant varieties and fat points.
  {Z}ero-dimensional schemes and applications - {N}aples, 2000.
\newblock 123:133--147, 2002.

\bibitem[CGG02b]{CGG2_Segre}
\bysame. 
\newblock Ranks of tensors, secant varieties of {S}egre varieties and fat
  points.
\newblock {\em Linear Algebra Appl.}, 355:263--285, 2002.

\bibitem[CGG03]{CGG2.5}
\bysame. 
\newblock Erratum to: ``{R}anks of tensors, secant varieties of {S}egre
  varieties and fat points'' [{L}inear {A}lgebra {A}ppl.\ {\bf 355} (2002),
  263--285; {MR}1930149 (2003g:14070)].
\newblock {\em Linear Algebra Appl.}, 367:347--348, 2003.

\bibitem[CGG05a]{CGG3_SegreVeronese}
\bysame. 
\newblock Higher secant varieties of {S}egre-{V}eronese varieties.
\newblock In {\em Projective varieties with unexpected properties}, pages
  81--107. Berlin: Walter de Gruyter GmbH \& Co. KG, 2005.

\bibitem[CGG05b]{CGG4_Segre}
\bysame. 
\newblock Higher secant varieties of the {S}egre varieties {$\Bbb P\sp
  1\times\dots\times\Bbb P\sp 1$}.
\newblock {\em J. Pure Appl. Algebra}, 201(1-3):367--380, 2005.

\bibitem[CGG05c]{CGG6_Grassman}
\bysame. 
\newblock Secant varieties of {G}rassmann varieties.
\newblock {\em Proc. Amer. Math. Soc.}, 133(3):633--642 (electronic), 2005.

\bibitem[CGG08]{CGG5_rational}
\bysame. 
\newblock On the ideals of secant varieties to certain rational varieties.
\newblock {\em J. Algebra}, 319(5):1913--1931, 2008.

\bibitem[CGG11]{CGG_P1s}
\bysame. 
\newblock Secant varieties of {$\Bbb P^1\times\dots\times\Bbb P^1$}
  ({$n$}-times) are not defective for {$n\geq 5$}.
\newblock {\em J. Algebraic Geom.}, 20(2):295--327, 2011.

\bibitem[CCD{\etalchar{+}}11]{CCDRS_Discriminants}
E.~{Cattani}, M.A. {Cueto}, A.~{Dickenstein}, S.~{Di Rocco}, and
  B.~{Sturmfels}.
\newblock {Mixed Discriminants}.
\newblock {\em Math. Z.} 274(3-4):761--778, 2013.

\bibitem[CR11]{Causa_rank_binary}
A.~Causa and R.~Re.
\newblock On the maximum rank of a real binary form.
\newblock {\em Annali di Matematica Pura ed Applicata}, 190:55--59, 2011.

\bibitem[CC02]{ChiCil_WeakDefective}
L.~Chiantini and C.~Ciliberto.
\newblock Weakly defective varieties.
\newblock {\em Trans. Amer. Math. Soc.}, 354(1):151--178 (electronic), 2002.

\bibitem[CC06]{ChiCil06_ksecant}
\bysame. 
\newblock On the concept of {$k$}-secant order of a variety.
\newblock {\em J. London Math. Soc. (2)}, 73(2):436--454, 2006.

\bibitem[Chi04]{Chipalkatti04a_Waring}
J.~Chipalkatti.
\newblock The {W}aring loci of ternary quartics.
\newblock {\em Experiment. Math.}, 13(1):93--101, 2004.

\bibitem[CGH{\etalchar{+}}05]{CGHRR}
C.~Ciliberto, A.~V. Geramita, B.~Harbourne, R.~M. Mir{\'o}-Roig, and
  K.~Ranestad, editors.
\newblock {\em Projective varieties with unexpected properties}. Berlin: Walter
  de Gruyter GmbH \& Co. KG, 2005.
\newblock A volume in memory of Giuseppe Veronese.

\bibitem[CS11b]{ComSei11_binary}
G.~Comas and M.~Seiguer.
\newblock On the rank of a binary form.
\newblock {\em Foundations of Computational Mathematics}, 11:65--78, 2011.

\bibitem[Com00]{Comon2_tensor}
P.~Comon.
\newblock Tensors decompositions. {S}tate of the art and applications.
\newblock {\em IMA Conference of Mathematics in Signal Processing}, December
  18-20, 2000.

\bibitem[Com04]{Comon1_canonical}
\bysame. 
\newblock Canonical tensor decompositions.
\newblock {\em Research Report RR-2004-17, I3S}, 17 June 2004.

\bibitem[CR06]{ComoR06_blind}
P.~Comon and M.~Rajih.
\newblock Blind identification of under-determined mixtures based on the
  characteristic function.
\newblock {\em Signal Processing}, 86(9):2271--2281, 2006.


\bibitem[CGLM08]{CGLM_Waring}
P.~Comon, G.~Golub, L.-H. Lim, and B.~Mourrain.
\newblock Symmetric tensors and symmetric tensor rank.
\newblock {\em SIAM J. Matrix Anal. Appl.}, 30(3):1254--1279, 2008.



\bibitem[CM96]{Comon-Mourrain_Waring}
P.~Comon and B.~Mourrain.
\newblock Decomposition of quantics in sums of powers of linear forms.
\newblock {\em Signal Processing}, 53(2):93--107, September 1996.
\newblock Special issue on High-Order Statistics.

\bibitem[CO12]{ComonOtt_Binary}
P. Comon and G. Ottaviani.
\newblock On the typical rank of real binary forms.
\newblock {\em Linear and Multilinear Algebra}, 60(6):657--667, 2012.


\bibitem[CtBDLC09]{ComonTenBerg_Typical}
P.~Comon, J.~M.~F. ten Berge, L.~De~Lathauwer, and J.~Castaing.
\newblock Generic and typical ranks of multi-way arrays.
\newblock {\em Linear Algebra Appl.}, 430(11-12):2997--3007, 2009.

\bibitem[CLO07]{CLO_text}
D.~Cox, J.~Little, and D.~O'Shea.
\newblock {\em Ideals, varieties, and algorithms: An introduction to
  computational algebraic geometry and commutative algebra}.
\newblock Undergraduate Texts in Mathematics. New York: Springer, third
  edition, 2007.

\bibitem[{CoC}]{CoCoA}
{CoCoA}Team.
\newblock {{\hbox{\rm C\kern-.13em o\kern-.07em C\kern-.13em o\kern-.15em A}}}:
  a system for doi\ ng {C}omputations in {C}ommutative {A}lgebra.
\newblock Available at \url{ http://cocoa.dima.unige.it}.



\bibitem[DBL07]{DBLP:conf/coco/2007}
{\em 22nd Annual IEEE Conference on Computational Complexity (CCC 2007), 13-16
  June 2007, San Diego, California, USA}. IEEE Computer Society, 2007.

\bibitem[DGPS10]{Singular}
W.~Decker, G.-M. Greuel, G.~Pfister, and H.~Sch{\"o}nemann.
\newblock {\sc Singular} {3-1-1} --- {A} computer algebra system for polynomial
  computations.
\newblock {\tt http://www.singular.uni-kl.de.}



\bibitem[DLC08]{LathCast_blind}
L.~De~Lathauwer and J.~Castaing.
\newblock Blind identification of underdetermined mixtures by simultaneous
  matrix diagonalization.
\newblock {\em Signal Processing, IEEE Transactions on}, 56(3):1096--1105,
  March 2008.

\bibitem[DLCC07]{LathCastCard_cumulant}
\bysame. 
\newblock Fourth-order cumulant-based blind identification of underdetermined
  mixtures.
\newblock {\em Signal Processing, IEEE Transactions on}, 55(6):2965--2973, June
  2007.

\bibitem[DLFDMV02]{DFDV02_block_decomp}
L.~De~Lathauwer, C.~F\'evotte, B.~De~Moor, and J.~Vandewalle.
\newblock Jacobi algorithm for joint block diagonalization in blind
  identification.
\newblock In {\em Proc. 23rd Symp. on Information Theory in the Benelux}, pages
  155--162, Louvain-la-Neuve, Belgium, May 2002.

\bibitem[dSL08]{deSilva_Lim_ill_posed}
V.~de~Silva and L.H. Lim.
\newblock Tensor rank and the ill-posedness of the best low-rank approximation
  problem.
\newblock {\em SIAM J. Matrix Anal. Appl.}, 30:1084--1127, September 2008.

\bibitem[Djo83]{Djokovic}
Dragomir~{\v{Z}}. Djokovi{\'c}.
\newblock Closures of equivalence classes of trivectors of an eight-dimensional
  complex vector space.
\newblock {\em Canad. Math. Bull.}, 26(1):92--100, 1983.

\bibitem[DK11]{draisma-kuttler_bounded}
J.~{Draisma} and J.~{Kuttler}.
\newblock {Bounded-rank tensors are defined in bounded degree}.
\newblock {\em ArXiv e-prints}, March 2011.

\bibitem[DSS07]{Drton-Sturmfels-Sullivant}
M.~Drton, B.~Sturmfels, and S.~Sullivant.
\newblock Algebraic factor analysis: tetrads, pentads and beyond.
\newblock {\em Probab. Theory Related Fields}, 138(3-4):463--493, 2007.

\bibitem[DSS09]{DMS_book}
\bysame. 
\newblock {\em {Lectures on algebraic statistics.}}
\newblock {Oberwolfach Seminars 39. Basel: Birkh\"auser. viii, 271~p.}, 2009.


\bibitem[Eis95]{Eisenbud_tome}
D.~Eisenbud.
\newblock {\em Commutative algebra with a view toward algebraic geometry.} 
\newblock Springer-Verlag, New York, 1995.
\newblock  Graduate Texts in
  Mathematics, No.150.

\bibitem[Eis05]{Eisenbud_syzygies}
\bysame. 
\newblock {\em The geometry of syzygies, second course in commutative algebra and algebraic geometry.} 
\newblock Springer-Verlag, New York, 2005. Graduate Texts in
  Mathematics, No. 229.
\newblock 

\bibitem[EC18]{Enriques_book}
F.~Enriques and O.~Chisini.
\newblock {\em Lezioni sulla teoria geometrica delle equazioni e delle funzioni
  algebriche}.
\newblock Zanichelli, Bologna, 1918.

\bibitem[ERSS05]{ERSS_phylo}
N. Eriksson, K. Ranestad, B. Sturmfels, and S. Sullivant.
\newblock Phylogenetic algebraic geometry.
\newblock In {\em Projective varieties with unexpected properties}, pages
  237--255. Walter de Gruyter GmbH \& Co. KG, Berlin, 2005.

\bibitem[EV10]{ErmanVelasco_smoothable}
D. Erman and M. Velasco.
\newblock A syzygetic approach to the smoothability of zero-dimensional
  schemes.
\newblock {\em Adv. in Math.}, 224:1143--1166, 2010.

\bibitem[FOV99]{FlennerOCarrol_Joins}
H.~Flenner, L.~O'Carroll, and W.~Vogel.
\newblock {\em Joins and intersections}.
\newblock Springer Monographs in Mathematics. Springer-Verlag, Berlin, 1999.





\bibitem[Fri12]{Friedland_Typical}
S. Friedland.
\newblock On the generic and typical ranks of 3-tensors.
\newblock {\em Linear Algebra Appl.}, 436(3):478--497, 2012.

\bibitem[Fri13]{Friedland2010_salmon}
\bysame. 
\newblock On tensors of border rank {$l$} in {$\Bbb{C}^{m\times n\times l}$}.
\newblock {\em Linear Algebra Appl.}, 438(2):713--737, 2013.

\bibitem[FG12]{Friedland-Gross2011_salmon}
S. Friedland and E. Gross.
\newblock A proof of the set-theoretic version of the salmon conjecture.
\newblock {\em J. Algebra}, 356:374--379, 2012.

\bibitem[FOS12]{FrobergOttavianiShapiro_Waring}
R. Fr\"oberg, G. Ottaviani, and B. Shapiro.
\newblock On the waring problem for polynomial rings.
\newblock {\em Proceedings of the National Academy of Sciences},
  109(15):5600--5602, 2012.

\bibitem[FH91]{FultonHarris}
W.~Fulton and J.~Harris.
\newblock {\em Representation Theory, a First Course.}
\newblock New York: Springer-Verlag, 1991.  Graduate Texts in
  Mathematics, No. 129.


\bibitem[Gan59]{gantmacher}
F.~R. Gantmacher.
\newblock {\em Applications of the theory of matrices}.
\newblock Translated by J. L. Brenner, with the assistance of D. W. Bushaw and
  S. Evanusa. Interscience Publishers, Inc., New York, 1959.

\bibitem[GSS05]{GSS}
L.~Garcia, M.~Stillman, and B.~Sturmfels.
\newblock Algebraic geometry of {B}ayesian networks.
\newblock {\em J. Symbolic Comput.}, 39(3-4):331--355, 2005.
\bibitem[GS]{M2}
D.~R. Grayson and M.~E. Stillman. Macaulay2, a software system for research in algebraic
  geometry.
\newblock Available at \url{http://www.math.uiuc.edu/Macaulay2/}.


\bibitem[GKZ92]{GKZ92}
I.~M. Gel{\cprime}fand, M.~M. Kapranov, and A.~V. Zelevinsky.
\newblock Hyperdeterminants.
\newblock {\em Adv. Math.}, 96(2):226--263, 1992.

\bibitem[GKZ94]{GKZ}
\bysame. 
\newblock {\em Discriminants, resultants, and multidimensional determinants}.
\newblock Mathematics: Theory \& Applications. Boston: Birkh\"auser, Boston,
  MA, 1994.
  
\bibitem[GZK89]{GKZ89}
I.~M. Gel{\cprime}fand, A.~V. Zelevinski{\u\i}, and M.~M. Kapranov.
\newblock Projective-dual varieties and hyperdeterminants.
\newblock {\em Dokl. Akad. Nauk SSSR}, 305(6):1294--1298, 1989.

\bibitem[Ger96]{MR1381732}
A.~V. Geramita.
\newblock Inverse systems of fat points: {W}aring's problem, secant varieties
  of {V}eronese varieties and parameter spaces for {G}orenstein ideals.
\newblock In {\em The {C}urves {S}eminar at {Q}ueen's, {V}ol.\ {X} ({K}ingston,
  {ON}, 1995)}, volume 102 of {\em Queen's Papers in Pure and Appl. Math.},
  pages 2--114. Queen's Univ., Kingston, ON, 1996.



\bibitem[GW98]{GoodWall}
R.~Goodman and N.~Wallach.
\newblock {\em Representations and invariants of the classical groups},
  volume~68 of {\em Encyclopedia of Mathematics and its Applications}.
\newblock Cambridge University Press, 1998.



\bibitem[Har92]{Harris}
J.~Harris.
\newblock {\em Algebraic geometry: A First Course}. 
\newblock New York: Springer-Verlag, 1992.
\newblock Graduate Texts in Mathematics, No. 133.


\bibitem[Har77]{Hartshorne_tome}
R. Hartshorne.
\newblock {\em Algebraic geometry}.
\newblock Springer-Verlag, New York, 1977.
\newblock Graduate Texts in Mathematics, No. 52.

\bibitem[HSW11]{regen}
J.D. Hauenstein, A.J. Sommese, and C.W. Wampler.
\newblock Regeneration homotopies for solving systems of polynomials.
\newblock {\em Math. Comp.}, 80:345--377, 2011.

\bibitem[HL09]{HillarLim09mosttensor}
C.~J. Hillar and L.H. Lim.
\newblock {Most tensor problems are NP hard}. {\em ArXiv e-prints}, November  2009.

\bibitem[HHLQ13]{Hu_Edet}
S. Hu, Z.H. Huang, C. Ling, and L. Qi.
\newblock On determinants and eigenvalue theory of tensors.
\newblock {\em Journal of Symbolic Computation}, 50(0):508--531, 2013.



\bibitem[HSYY08]{HSYY_hyperdet}
P.~Huggins, B.~Sturmfels, J.~Yu, and D.~Yuster.
\newblock The hyperdeterminant and triangulations of the 4-cube.
\newblock {\em Math. Comp.}, 77(263):1653--1679, 2008.

\bibitem[IK99]{IarrobinoKanev_text}
A.~Iarrobino and V.~Kanev.
\newblock {\em Power sums, {G}orenstein algebras, and determinantal loci},
  volume 1721 of {\em Lecture Notes in Mathematics}.
\newblock Springer-Verlag, Berlin, 1999.
\newblock Appendix C by Iarrobino and Steven L. Kleiman.

\bibitem[Kan99]{Kanev_catalecticant}
V~Kanev.
\newblock Chordal varieties of {V}eronese varieties and catalecticant matrices.
\newblock {\em J. Math. Sci. (New York)}, 94(1):1114--1125, 1999.
\newblock Algebraic geometry, 9.

\bibitem[KM11]{KoldaMayo_tensor_eigenpairs}
T.~G. Kolda and J.~R. Mayo.
\newblock Shifted power method for computing tensor eigenpairs.
\newblock {\em SIAM J. Matrix Anal. Appl.}, 32(4):1095--1124, 2011.

\bibitem[Lan06]{Lan_2by2_matrix}
J.~M. Landsberg.
\newblock The border rank of the multiplication of {$2\times2$} matrices is
  seven.
\newblock {\em J. Amer. Math. Soc.}, 19(2):447--459 (electronic), 2006.

\bibitem[Lan08]{Lan_matrix_mult}
\bysame. 
\newblock Geometry and the complexity of matrix multiplication.
\newblock {\em Bull. Amer. Math. Soc. (N.S.)}, 45(2):247--284, 2008.


\bibitem[{Lan}12a]{Landsberg_explicit}
\bysame. 
\newblock {Explicit tensors of border rank at least 2n-1}.
\newblock {\em ArXiv e-prints}, September 2012.

\bibitem[{Lan}12b]{Landsberg_NewLower}
\bysame. 
\newblock {New lower bounds for the rank of matrix multiplication}.
\newblock {\em ArXiv e-prints}, June 2012.

\bibitem[Lan12c]{LandsbergTensorBook}
\bysame. 
\newblock {\em Tensors: geometry and applications}, volume 128 of {\em Graduate
  Studies in Mathematics}.
\newblock American Mathematical Society, Providence, RI, 2012.

\bibitem[LM04]{LanMan04_Seg}
J.~M. Landsberg and L.~Manivel.
\newblock On the ideals of secant varieties of {S}egre varieties.
\newblock {\em Found. Comput. Math.}, 4(4):397--422, 2004.

\bibitem[LM08]{LanMan08_Strassen}
\bysame. 
\newblock  Generalizations of {S}trassen's equations for secant varieties of
  {S}egre varieties.
\newblock {\em Comm. Algebra}, 36(2):405--422, 2008.

\bibitem[LO11a]{LanOtt_NewLower}
J.~M. {Landsberg} and G.~{Ottaviani}.
\newblock {New lower bounds for the border rank of matrix multiplication}.
\newblock {\em ArXiv e-prints}, December 2011.

\bibitem[LO11b]{LanOtt11_Equations}
\bysame. 
\newblock Equations for secant varieties of veronese and other varieties.
\newblock {\em Annali di Matematica Pura ed Applicata}, pages 1--38, 2011.

\bibitem[LT10]{LandsbergTeitler_SymRank}
J.~M. Landsberg and Z. Teitler.
\newblock On the ranks and border ranks of symmetric tensors.
\newblock {\em Found. Comput. Math.}, 10(3):339--366, 2010.

\bibitem[LW07a]{LanWey_Tan}
J.~M. Landsberg and J. Weyman.
\newblock On tangential varieties of rational homogeneous varieties.
\newblock {\em J. Lond. Math. Soc. (2)}, 76(2):513--530, 2007.

\bibitem[LW07b]{LanWey_Seg}
\bysame. 
\newblock On the ideals and singularities of secant varieties of {S}egre
  varieties.
\newblock {\em Bull. Lond. Math. Soc.}, 39(4):685--697, 2007.

\bibitem[LW09]{LanWey_secant_CHSS}
\bysame. 
\newblock On secant varieties of compact {H}ermitian symmetric spaces.
\newblock {\em J. Pure Appl. Algebra}, 213(11):2075--2086, 2009.


\bibitem[Lat06]{lathauwer:642_matrix_diag}
L.~De Lathauwer.
\newblock A link between the canonical decomposition in multilinear algebra and
  simultaneous matrix diagonalization.
\newblock {\em SIAM Journal on Matrix Analysis and Applications},
  28(3):642--666, 2006.

\bibitem[LdB08]{DD08_blind}
L.~De Lathauwer and A.~de~Baynast.
\newblock Blind deconvolution of {DS-CDMA} signals by means of decomposition in
  rank-{$(1,L,L)$} terms.
\newblock {\em {IEEE} Trans. Signal Processing}, 56(4):1562--1571, 2008.

\bibitem[LMV04]{lathauwer_Schur}
L.~De Lathauwer, B.~De Moor, and Joos Vandewalle.
\newblock Computation of the canonical decomposition by means of a simultaneous
  generalized {S}chur decomposition.
\newblock {\em SIAM Journal on Matrix Analysis and Applications},
  26(2):295--327, 2004.

\bibitem[LLT]{HOM4PS}
T.L. Lee, T.Y. Li, and C.H. Tsai.
\newblock  
  HOM4PS-2.0: 
  A software package for solving polynomials systems by
  the polyhedral homotopy continuation method. Available at \url{http://www.mth.msu.edu/~li/}, 2010. 

\bibitem[Li03]{Li_homotopy}
T.Y. Li.
\newblock Numerical solution of polynomial systems by homotopy continuation
  methods.
\newblock volume~XI of {\em Handbook of Numerical Analysis, Special Volume:
  Foundations of Computational Mathematics}, pages 209--304. North-Holland,
  2003.

\bibitem[Lim05]{Lim05_evectors}
L.H. Lim.
\newblock Singular values and eigenvalues of tensors: a variational approach.
\newblock In {\em Computational Advances in Multi-Sensor Adaptive Processing,
  2005 1st IEEE International Workshop on}, pages 129--132, December 2005.





\bibitem[Mac1927]{Mac1927}
F.~S. Macaulay. 
\newblock {S}ome {P}roperties of {E}numeration in the {T}heory of {M}odular {S}ystems.
\newblock {\em Proc. London Math. Soc.},  S2-26(1):531--555, 1927.


\bibitem[Mel06]{Mella06_Waring}
M.~Mella.
\newblock Singularities of linear systems and the {W}aring problem.
\newblock {\em Trans. Amer. Math. Soc.}, 358(12):5523--5538 (electronic), 2006.

\bibitem[Mel09]{Mella09_Waring}
\bysame. 
\newblock Base loci of linear systems and the {W}aring problem.
\newblock {\em Proc. Amer. Math. Soc.}, 137(1):91--98, 2009.


\bibitem[NQWW07]{NiQiWang_Epoly}
G.~Ni, L.~Qi, F~Wang, and Y.~Wang.
\newblock The degree of the {E}-characteristic polynomial of an even order
  tensor.
\newblock {\em J. Math. Anal. Appl.}, 329(2):1218--1229, 2007.

\bibitem[Oed11]{OedingTan}
L.~Oeding.
\newblock Set-theoretic defining equations of the tangential variety of the
  Segre variety.
\newblock {\em Journal of Pure and Applied Algebra}, 215(6):1516 -- 1527, 2011.

\bibitem[Oed12]{Oeding_Hyperdet}
\bysame. 
\newblock Hyperdeterminants of polynomials.
\newblock {\em Adv. Math.}, 231(3-4):1308--1326, 2012.

\bibitem[OO13]{OedOtt13_Waring}
L. Oeding and G. Ottaviani.
\newblock Eigenvectors of tensors and algorithms for waring decomposition.
\newblock {\em Journal of Symbolic Computation}, (54):9--35, 2013.

\bibitem[OSS80]{OSS_text}
C.~Okonek, M.~Schneider, and H.~Spindler.
\newblock {\em Vector bundles on complex projective spaces}, volume~3 of {\em
  Progress in Mathematics}.
\newblock Birkh\"auser Boston, Mass., 1980.

\bibitem[OS13]{OttSturm_Eigenvectors}
G.~{Ottaviani} and B.~{Sturmfels}.
\newblock {Matrices with Eigenvectors in a Given Subspace}.
\newblock {\em Proc. of the American Math. Soc.}, 141:1219--1232, 2013.


\bibitem[Ott07]{Ott07_Luroth}
G.~Ottaviani.
\newblock Symplectic bundles on the plane, secant varieties and {L}\"uroth
  quartics revisited.
\newblock In {\em Vector bundles and low codimensional subvarieties: state of
  the art and recent developments}, volume~21 of {\em Quad. Mat.}, pages
  315--352. Dept. Math., Seconda Univ. Napoli, Caserta, 2007.

\bibitem[Ott09]{Ottaviani09_Waring}
\bysame. 
\newblock An invariant regarding {W}aring's problem for cubic polynomials.
\newblock {\em Nagoya Math. J.}, 193:95--110, 2009.

\bibitem[PS05]{Pachter-Sturmfels}
L.~Pachter and B.~Sturmfels, editors.
\newblock {\em Algebraic statistics for computational biology}.
\newblock New York: Cambridge University Press, 2005.

\bibitem[Qi05]{Qi05_eigen}
L.~Qi.
\newblock Eigenvalues of a real supersymmetric tensor.
\newblock {\em J. Symbolic Comput.}, 40(6):1302--1324, 2005.

\bibitem[{Rai}10]{RaicuCat}
C.~{Raicu}.
\newblock { $3 \times 3 $ Minors of Catalecticants}.
\newblock {\em ArXiv e-prints}, November 2010.

\bibitem[Rai11]{Raicu_thesis}
\bysame. 
\newblock {\em Secant {V}arieties of {S}egre-{V}eronese {V}arieties}.
\newblock ProQuest LLC, Ann Arbor, MI, 2011.
\newblock Thesis (Ph.D.), University of California, Berkeley.

\bibitem[{Rai}12]{RaicuGSS}
\bysame. 
\newblock Secant varieties of Segre--Veronese varieties.
\newblock {\em Algebra and Number Theory}, 6-8:1817--1868, 2012.

\bibitem[RS00]{RanestadSchreyer_VSP}
K.~Ranestad and F.O. Schreyer.
\newblock Varieties of sums of powers.
\newblock {\em J. Reine Angew. Math.}, 525:147--181, 2000.

\bibitem[SS09]{SidmanSullivant_Prolongations}
J. Sidman and S. Sullivant.
\newblock Prolongations and computational algebra.
\newblock {\em Canad. J. Math.}, 61(4):930--949, 2009.

\bibitem[SV11]{SidmanVermeire_Secant}
J. Sidman and P. Vermeire.
\newblock Equations defining secant varieties: geometry and computation.
\newblock In {\em Combinatorial aspects of commutative algebra and algebraic
  geometry}, volume~6 of {\em Abel Symp.}, pages 155--174. Springer, Berlin,
  2011.

\bibitem[SW05]{SomeseWampler05}
A.J. Sommese and C.W. Wampler.
\newblock {\em Numerical solution of polynomial systems arising in engineering
  and science}.
\newblock World Scientific, Singapore, 2005.

\bibitem[SSM13]{SumiSakata_Typical}
T. Sumi, T. Sakata, and M. Miyazaki.
\newblock Typical ranks for {$m\times n\times (m-1)n$} tensors with {$m\leq
  n$}.
\newblock {\em Linear Algebra Appl.}, 438(2):953--958, 2013.

\bibitem[Str83]{Strassen83_rank}
V.~Strassen.
\newblock Rank and optimal computation of generic tensors.
\newblock {\em Linear Algebra Appl.}, (52-53):645--685, 1983.

\bibitem[Stu09]{SturmfelsOpenProblems}
B.~Sturmfels.
\newblock Open problems in algebraic statistics.
\newblock In S.~Sullivant and B.~Sturmfels, editors, {\em Emerging Applications
  of Algebraic Geometry}, volume 149 of {\em The IMA Volumes in Mathematics and
  its Applications}, pages 1--13. Springer New York, 2009.



\bibitem[Ter11]{Terracini}
A.~Terracini.
\newblock Sulle vk per cui la variet`a degli sh h + 1-secanti ha dimensione
  minore dell'ordinario.
\newblock {\em Rend. Circ. Mat. Palermo}, Selecta vol. I:392--396, 1911.

\bibitem[Tre12]{SIAM_news}
N. Trefethen.
\newblock The smart money is on numerical analysts.
\newblock 45(9), 2012.

\bibitem[vLCL92]{lie}
M.~A.~A. van Leeuwen, A.~M. Coehn, and B.~Lisser.
\newblock {\em {LiE}, A Package for {Lie} Group Computations}.
\newblock Computer Algebra Nederland, 1992.


\bibitem[Ver]{PHC}
J.~Verschelde.
\newblock {PHCpack}: a general-purpose solver for
  polynomial systems by homotopy continuation.  Available at \url{http://www.math.uic.edu/~jan/PHCpack/phcpack.html}, 2010.

\bibitem[V{\`E}78]{VinbergElashvili_classification}
{\`E}.~B. Vinberg and A.~G. {\`E}la{\v{s}}vili.
\newblock A classification of the three-vectors of nine-dimensional space.
\newblock {\em Trudy Sem. Vektor. Tenzor. Anal.}, 18:197--233, 1978.


\bibitem[Wey97]{Weyl}
H. Weyl.
\newblock {\em The classical groups}.
\newblock Princeton Landmarks in Mathematics. Princeton University Press,
  Princeton, NJ, 1997.
\newblock Their invariants and representations, Fifteenth printing, Princeton
  Paperbacks.

\bibitem[Wey03]{Weyman}
J.~Weyman.
\newblock {\em Cohomology of vector bundles and syzygies}, volume 149 of {\em
  Cambridge Tracts in Mathematics}.
\newblock Cambridge University Press, 2003.

\bibitem[WZ94]{WeymanZelevinsky_mult}
J.~Weyman and A.~Zelevinsky.
\newblock Multiplicative properties of projectively dual varieties.
\newblock {\em Manuscripta Math.}, 82(2):139--148, 1994.

\bibitem[WZ96]{WeymanZelevinsky_sing}
\bysame. 
\newblock Singularities of hyperdeterminants.
\newblock {\em Ann. Inst. Fourier (Grenoble)}, 46(3):591--644, 1996.

\bibitem[Wil12]{Williams12}
V.V. Williams.
\newblock Multiplying matrices faster than coppersmith-winograd.
\newblock In Howard~J. Karloff and Toniann Pitassi, editors, {\em STOC}, pages
  887--898. ACM, 2012.

\bibitem[Zak93]{Zak}
F.~L. Zak.
\newblock {\em Tangents and secants of algebraic varieties}, volume 127 of {\em
  Translations of Mathematical Monographs}.
\newblock Providence: American Mathematical Society, 1993.
\newblock Translated from the Russian manuscript by the author.

\end{thebibliography}

\newcommand{\etalchar}[1]{$^{#1}$}
\def\Dbar{\leavevmode\lower.6ex\hbox to 0pt{\hskip-.23ex \accent"16\hss}D}
  \def\cprime{$'$}
  \providecommand{\bysame}{\leavevmode\hbox to3em{\hrulefill}\thinspace}

\end{document}